\documentclass[12pt]{amsart}

\usepackage{pinlabel}
\usepackage{amsmath,amsfonts,amssymb,latexsym,mathrsfs,stmaryrd}
\usepackage{color}
\usepackage[margin=3cm]{geometry}
\usepackage{graphicx}
\usepackage{graphics}
\usepackage{enumerate}
\usepackage{amscd} 
\usepackage{amsxtra}
\usepackage{epic,eepic}
\usepackage{hyperref} 
\usepackage{enumitem}
\usepackage{MnSymbol}
\usepackage{dsfont}
\usepackage[all]{xy}
\usepackage[alphabetic]{amsrefs}

\newtheorem{theorem}{Theorem}[section]

\newtheorem{proposition}[theorem]{Proposition}

\newtheorem{thm}[theorem]{Theorem}

\theoremstyle{definition} 
\newtheorem{defn}[theorem]{Definition}
\newtheorem{definition}[theorem]{Definition}
\newtheorem{example}[theorem]{Example}

\newtheorem{remark}[theorem]{Remark}
\newtheorem{rmk}[theorem]{Remark}

\newcommand{\on}{\operatorname} 
\newcommand{\HH}{\mathfrak H}

\title[DHT for toric complete intersections]{The Doran-Harder-Thompson Conjecture for toric complete intersections}

\author{Charles F. Doran}
\address{Department of Mathematical and Statistical Sciences, University of Alberta and Center of Mathematical Sciences and Applications, Harvard University}
\email{charles.doran@ualberta.ca}

\author{Jordan Kostiuk}
\address{Department of Mathematics, Brown University}
\email{jordan\_kostiuk@brown.edu}

\author{Fenglong You}
\address{Department of Mathematical and Statistical Sciences, University of Alberta and Fields Institute for Research in the Mathematical Sciences}
\email{fenglong@ualberta.ca}

\begin{document}

\maketitle

\begin{abstract}
    Given a Tyurin degeneration of a Calabi-Yau complete intersection in a toric variety, we prove gluing formulas relating the generalized functional invariants, periods, and $I$-functions of the mirror Calabi-Yau family and those of the two mirror  Landau-Ginzburg models. Our proof makes explicit the ``gluing/splitting'' of fibrations in the Doran-Harder-Thompson  mirror conjecture.  Our gluing formula implies an identity, obtained by composition with their respective mirror maps, that relates the absolute Gromov-Witten invariants for the Calabi-Yaus and relative Gromov-Witten invariants for the quasi-Fanos.  
\end{abstract}
\tableofcontents

\section{Introduction}

\subsection{The Doran-Harder-Thompson conjecture}

Classical mirror symmetry is a conjecture relating properties of a Calabi-Yau variety and properties of its mirror Calabi-Yau variety. The duality has been generalized to Fano varieties. By \cite{EHX}, the mirror to a Fano variety $X$ is a Landau-Ginzburg model $(X^\vee,W)$ such that $X^\vee$ is a K\"ahler manifold and $W:X^\vee\rightarrow \mathbb C$ is a proper map called the superpotential. 

Mirror symmetry for Landau-Ginzburg models can be generalized to varieties beyond the Fano case. In particular, the Landau-Ginzburg model for a quasi-Fano variety is defined by A.~Harder \cite{Harder16}. A smooth variety $X$ is quasi-Fano if its anticanonical linear system contains a smooth Calabi-Yau member and $H^i(X,\mathcal O_X)=0$ for all $i>0$. A degeneration of Calabi-Yau varieties is given by $\mathcal V\rightarrow \Delta$, where $\Delta\subset \mathbb C$ is the unit disk. The degeneration is called a Tyurin degeneration if the total space $\mathcal V$ is smooth and the central fiber consists of two quasi-Fano varieties $X_1$ and $X_2$ which meet normally along a common smooth anticanonical (Calabi-Yau) divisor $X_0$.

Motivated in part by a question of A.~Tyurin \cite{Tyurin04}, Doran-Harder-Thompson formulated a conjecture relating mirror symmetry for a Calabi-Yau variety $X$ arising as a smooth fiber of a Tyurin degeneration $\mathcal V \rightarrow \Delta$, and mirror symmetry for the quasi-Fano varieties $X_1$ and $X_2$ in the central fiber of $\mathcal V$. The Doran-Harder-Thompson conjecture states that one should be able to glue the Landau-Ginzburg models $W_i:X^\vee_i \rightarrow \mathbb C$ of the pair $(X_i,X_0)$ for $i=1,2$ to a Calabi-Yau variety $X^\vee$ which is mirror to $X$ with the fibers of the superpotentials gluing to a fibration $W:X^\vee\rightarrow \mathbb P^1$. Furthermore, the compact fibers of the Landau-Ginzburg models consist of Calabi-Yau manifolds mirror to the common anticanonical divisor $X_0$.  

An enormous amount of evidence has already been collected in favor of this conjecture.  The topological version of the conjecture has been proved in \cite{DHT} by computing the corresponding Euler numbers.  In this case the gluing could be viewed as taking place in the symplectic category, with monodromies aligning with mirror autoequivalences in the bounded derived category on the Tyurin degeneration side. In the case of elliptic curves, the conjecture is proved by Kanazawa \cite{Kanazawa17} via SYZ mirror symmetry. The original paper \cite{DHT} included evidence of compatibility of the DHT conjecture with the Dolgachev-Nikulin-Pinkham formulation of mirror symmetry for lattice polarized $K3$ surfaces.  This was refined by Doran-Thompson, who obtained mirror notions of lattice polarization for rational elliptic surfaces and for del Pezzo surfaces coming from, respectively, splitting of fibrations and Tyurin degeneration of lattice polarized $K3$ surfaces
\cite{doranMirrorSymmetryLattice2018}.
\subsection{The gluing formula}

In this paper, we study the Tyurin degeneration of complete intersections in toric varieties, and derive a gluing formula (in the complex category) for Landau-Ginzburg models. 

We first analyze two special examples, namely, a conifold transition $\tilde {Q}_5$ of the quintic threefold and the quintic threefold $Q_5$ itself. In these two examples, the mirror Calabi-Yau threefold is fibered by $M_2$-polarized (``mirror quartic'') $K3$ surfaces and the mirror LG models are $M_2$-polarized families of $K3$ surfaces. The key ingredients of the construction are the generalized functional invariants which determine the families of $M_2$-polarized $K3$ surfaces up to isomorphism. Therefore, gluing LG models must reduce to the gluing of generalized functional invariants. 

We show that the gluing formula for generalized functional invariants is simply a product relation between them. Furthermore, we use generalized functional invariants to define the holomorphic periods of LG models as pullbacks of the holomorphic period of  mirror quartic $K3$ surfaces. As shown by Doran-Malmendier \cite{doran_calabi-yau_2015} the holomorphic $3$-form periods of $K3$ surface fibered Calabi-Yau threefolds can also be computed from the generalized functional invariants and the holomorphic period of the fiber mirror quartic $K3$ surfaces.  They are expressed, equivalently, in terms of the Euler integral transform, the Hadamard product of power series, and the middle convolution of ordinary differential equations. As a result, the gluing formula for generalized functional invariants implies the gluing formula for periods as well as the gluing formula for the bases of solutions to the corresponding Picard-Fuchs equations (the full $I$-functions).

We generalize the gluing formula to toric complete intersections. After suitable identification of variables, we obtain the following relation among periods.

\begin{theorem}[= Theorem \ref{thm:gluing-period-general}]
Let $X$ be a Calabi-Yau complete intersection in a toric variety which corresponds to a nef partition. Consider the Tyurin degeneration of $X$ via a refinement of the nef partition. Let $X_1$ and $\tilde X_2$ be the corresponding quasi-Fano varieties, where $X_1$ and $X_2$ are toric complete intersections coming from the refinement of the nef partition and we blow up $X_2$ to obtain $\tilde X_2$. Let $X_0$ be the smooth anticanonical divisor which lies in the intersection of $X_1$ and $\tilde X_2$. The following Hadamard product relation holds
\[
f_0^{X}(q)\star_q f_0^{X_0}(q)=\frac{1}{2\pi i}\oint f_0^{X_1}(q,y)\star_q f_0^{\tilde{X}_2}(q,y)\frac{dy}{y}.
\]
\end{theorem}
We also have the Hadamard product relation among the solutions for the Picard-Fuchs equations. 
\begin{theorem}[= Theorem \ref{thm:general-case-I-function}]
The Hadamard product relation among the bases of solutions to Picard-Fuchs equations is
\[
I^{X}(q)\star_q I^{X_0}(q)=\frac{1}{2\pi i}\oint I^{X_1}(q,y)\star_q I^{\tilde{X}_2}(q,y)\frac{dy}{y}.
\]
We refer to Section \ref{sec:general-case} for explanation of the notation.
\end{theorem}

In this paper, we focus on toric complete intersections, but the gluing formulas work for more general Calabi-Yau manifolds. In particular, in a forthcoming paper, we will prove the gluing formulas for the Tyurin degenerations of Calabi-Yau threefolds mirror to the $M_N$-polarized $K3$ surface fibered Calabi-Yau threefolds classified and constructed in full in \cite{doran_calabi-yau_2017}.

\subsection{Gromov-Witten invariants}

According to the mirror theorem proved by Givental \cite{Givental98} and Lian-Liu-Yau \cite{LLY}, Gromov-Witten invariants of a Calabi-Yau variety are related to periods of the mirror. Following Givental \cite{Givental98}, on the A-side, we consider a generating function of Gromov-Witten invariants of the Calabi-Yau variety $X$, called the $J$-function. On the B-side, we consider the $I$-function, which encodes a basis of the corresponding Picard-Fuchs equation. The $J$-function and $I$-function are related by the mirror map. A mirror theorem for smooth pairs has recently been proved by Fan-Tseng-You \cite{FTY} using Givental's formalism for relative Gromov-Witten theory developed by Fan-Wu-You \cite{FWY}. Given a quasi-Fano variety $X$ and its anticanonical divisor $D$, similar to the absolute case, we can consider the relative $J$-function for the pair $(X,D)$. Under suitable assumptions, the relative $I$-function is related to the relative $J$-function via a relative mirror map. In general, the relative $I$-function lies in Givental's Lagrangian cone for relative Gromov-Witten theory defined in \cite{FWY}. 

The gluing formula for periods implies a relation among absolute Gromov-Witten invariants of the Calabi-Yau variety $X$ and the relative Gromov-Witten invariants of the pairs $(X_1,D)$ and $(X_2,D)$ via mirror maps. On the other hand, absolute Gromov-Witten invariants of $X$ and relative Gromov-Witten invariants of $(X_1,D)$ and $(X_2,D)$ are related by the degeneration formula. While we may consider our gluing formula as the B-model counterpart of the degeneration formula, the precise compatibility between the gluing formula and the degeneration formula is not yet known.

\subsection{Acknowledgment}
We would like to thank Andrew Harder, Hiroshi Iritani, Bumsig Kim, Melissa Liu, Yongbin Ruan, Alan Thompson, and Hsian-Hua Tseng for helpful discussions during various stages of this project.
C.~F.~D. acknowledges the support of the National Science and Engineering Research Council of Canada (NSERC).  F.~Y. is supported by a postdoctoral fellowship of NSERC and the Department of Mathematical and Statistical Sciences at the University of Alberta and a postdoctoral fellowship for the Thematic Program on Homological Algebra of Mirror Symmetry at the Fields Institute for Research in Mathematical Sciences.  The authors thank the Center of Mathematical Sciences and Applications (CMSA) at Harvard University where this work was completed and presented.

\section{Preparation}

\subsection{Generalized functional invariants}
In this section, we recall the definition of generalized functional invariants for threefolds. We also give a definition for generalized functional invariants in all dimensions.

We begin by reviewing some aspects of Kodaira's theory of ellipic surfaces, for which we refer to Kodaira's original papers \cites{kodaira_compact_1960,kodaira_compact_1963} for a complete treatment. 
Consider the following differential equation: 
\begin{equation}\label{Picard_Fuchs_J=t}
\frac{d^2f}{dt^2}+\frac{1}{t}\frac{df}{dt}+\frac{\frac{31}{144}t-\frac{1}{36}}{t^2(t-1)^2}f=0.
\end{equation}
Then, \eqref{Picard_Fuchs_J=t} is a Fuchsian differential equation with regular singularities at $t=0,1,\infty$. 
It admits a basis of solution $\omega_1,\omega_2$ for which the quotient $\tau=\frac{\omega_1}{\omega_2}$ defines a multi-valued function to the upper half-plane:
$$\mathbf{P}^1_t\Rrightarrow\mathfrak{H}.$$
The single-valued inverse function $J\colon\mathfrak{H}\to\mathbf{P}^1_t$ is the classical modular $J$-function.

The differential equation \eqref{Picard_Fuchs_J=t} is the Picard-Fuchs differential equation for the following family of elliptic curves: 
\begin{equation}\label{J=t Form}
y^2=4x^3-\frac{27t}{t-1}x-\frac{27t}{t-1}.	
\end{equation}
For each $t$, the elliptic curve $E_t$ is the elliptic curve with $J$-invariant equal to $t$. 

Given an arbitrary family of elliptic curves $f\colon\mathcal{E}\to T$, the period map $T\Rrightarrow\mathfrak{H}$ is the multi-valued function $\tau=\frac{\omega_1}{\omega_2}$ determined by choosing a suitable basis of period functions. 
The composition of the period map with the modular $J$-function is a \emph{rational function} $\mathcal{J}\colon T\to\mathbf{P}^1$ and is called \emph{Kodaira's functional invariant} associated to the family $f\colon\mathcal{E}\to T$. For each $t\in T$, $\mathcal{J}(t)$ is the $J$-invariant of the elliptic curve $E_t$. 

The \emph{homological invariant} of the family $f\colon\mathcal{E}\to T$ is the local system $\mathcal{G}=R^1f_*\mathbf{Z}$ of first cohomology groups of the family. 
Kodaira showed in \cite{kodaira_compact_1960} that the functional and homological invariants determine the isomorphism class of the elliptic surface. 
The functional invariant $\mathcal{J}$ is sufficient to determine the \emph{projective monodromy} representation of the homological invariant. 
It follows that the functional invariant, together with a representation $\pi_1(T)\to\{\pm 1\}$ determine the elliptic surface up to isomorphism.

A similar story holds for $M_N$-polarized families of $K3$ surfaces, as was proved in \cite{doran_calabi-yau_2017}. 
The period domain for such families is equal to $\mathfrak{H}$, and the moduli space is $\mathcal{M}_{M_N}\cong X_0(N)^+$, the quotient of the modular curve $X_0(N)$ by the Fricke involution. 
To any family $f\colon\mathcal{X}\to B$ of $M_N$-polarized $K3$ surface, we associate the \emph{generalized functional invariant} $g\colon B\to X_0(N)^+$ which is defined by sending each point $b\in B$ to the corresponding point in moduli of the fiber $f^{-1}(b)$. 
The generalized functional invariants associated to families of $M_N$-polarized $K3$ surfaces have even tighter control over the families: a family $f\colon\mathcal{X}\to B$ is determined up to isomorphism by the generalized functional invariant \cite{doran_calabi-yau_2017}. 

\begin{example}
Choose a coordinate $\lambda$ on the modular curve $X_0(2)^+$ for which $\lambda=0$ is the cusp, $\lambda=\frac{1}{256}$ is the order $2$ orbifold point and $\lambda=\infty$ is the order $4$ orbifold point. 
Let $\mathcal{X}_2$ denote the quartic mirror family of $K3$ surfaces defined by
$$\left\{\tilde{x}_1\tilde{x}_2\tilde{x}_3(\tilde{x}_1+\tilde{x}_2+\tilde{x}_3-1)+\lambda=0\right\}.$$

Then, the fibers of $\mathcal{X}_2$ are $M_2$-polarized $K3$ surfaces for $$\lambda\in X_0(2)^+-\left\{0,\frac{1}{256},\infty\right\}.$$ 
The singular fiber types of $\mathcal{X}_2$ are determined in \cite{doran_calabiyau_2016}. 
The singular fiber at the cusp $\lambda=0$ is a singular $K3$ surface of type $\textrm{III}$ containing $4$ components; the monodromy transformation of the corresponding variation of Hodge structure is maximally unipotent. 
At $\lambda=\frac{1}{256}$, the singular fiber is a singular $K3$ surface containing an $A_1$ singularity; the monodromy transformation on the VHS is conjugate to 
$$\begin{bmatrix}1&0&0\\
0&1&0\\
0&0&-1\end{bmatrix}.$$
The singular fiber at $\lambda=\infty$ is a singular fiber with $31$ components and the monodromy is conjugate to 
$$\begin{bmatrix}
i&0&0\\
0&-i&0\\
0&0&1
\end{bmatrix}.$$

According to the results described above, any $M_2$-polarized family of $K3$ surfaces is birational to the pull-back of $\mathcal{X}_2$ via the generalized functional invariant map $g\colon B\to X_0(2)^+$. 
The types of singular fibers appearing in the pull-back is determined by the ramification profile and is worked out in detail in \cite{doran_calabiyau_2016}.

A Picard-Fuchs operator corresponding to the VHS of the quartic mirror family is 
\begin{equation}
    \delta^3-256\lambda\left(\delta+\frac{1}{4}\right)\left(\delta+\frac{1}{2}\right)\left(\delta+\frac{3}{4}\right),
\end{equation}
where $\delta=t\frac{d}{dt}$. The holomorphic and logarithmic solutions to this ODE are given by
\begin{eqnarray*}
f_0^{\mathcal{X}_2}&=&\sum_{n=0}^\infty\frac{(4n)!}{(n!)^4}\lambda^n,\\
f_1^{\mathcal{X}_2}&=&f_0^{\mathcal{X}_2}\cdot\log{\lambda}+4\sum_{n=1}^{\infty}\frac{(4n)!}{(n!)^4}\cdot\left(\sum_{j=n+1}^{4n}\frac{1}{j}\right)\lambda^n.
\end{eqnarray*}

\end{example}

With this as motivation, we make the following definition:

\begin{definition}
Let $f\colon\mathcal{X}\to B$ be a family of Calabi-Yau manifolds. 
The morphism $g\colon B\to\mathcal{M}$, taking a point in the base to the corresponding point in the moduli space $\mathcal{M}$ of the fiber is called the \emph{generalized functional invariant} of the family. 
\end{definition}

\begin{remark}
The construction of $\mathcal{M}$ is strongly dependant on the type of Calabi-Yau fibers under consideration. 
If any confusion is likely to arise, this will be addressed on a case-by-case basis. The generalized functional invariants in this paper are rational functions, computed explicitly through toric mirror symmetry, coincides with the generalized functional invariants in \cite{doran_calabi-yau_2017} for $M_n$-polarized $K3$ surfaces. For general toric complete intersections, these rational functions are our generalization of the generalized functional invariants in \cite{doran_calabi-yau_2017}. 
\end{remark}

\begin{remark}
In this paper, we will be using Hori-Vafa mirrors because they can be used to describe both the mirrors of quasi-Fano varieties and the mirrors of Calabi-Yau complete intersections in toric varieties. The difference between generalized functional invariants in our paper and generalized functional invariants in previous literature (e.g. \cite{CDLNT}) is the following. First of all, in \cite{CDLNT}, authors considered Batyrev mirrors of Calabi-Yau complete intersections in toric varieties. Batyrev mirrors and Hori-Vafa mirrors are equivalent under change of variables. Secondly, generalized functional invariants in \cite{CDLNT}*{Section 4.1} are defined using the $M$-polarized K3 surface in \cite{CDLNT}*{Theorem 4.3}. While in this paper, we will consider generalized functional invariant maps to the mirror family of the Calabi-Yau $(n-1)$-fold in the Tyurin degeneration of Calabi-Yau $n$-fold. This is because we want to find the relation among generalized functional invariants of two LG models and mirror Calabi-Yau $n$-fold with respect to the mirror of the intersection of two quasi-Fanos. Therefore, after change of variables, generalized functional invariants in our paper will be the same as generalized functional invariants in \cite{CDLNT}.
\end{remark}

\subsection{Periods}\label{sec:period}
In this section, we define periods of LG models using generalized functional invariants. 
We first recall the definition of Landau-Ginzburg model for quasi-Fano varieties.

\begin{defn}[\cite{DHT}, Definition 2.1]
A Landau-Ginzburg model of a quasi-Fano variety $X$ is a pair $(X^\vee,W)$ consisting of a K\"ahler manifold $X^\vee$ satisfying $h^1(X^\vee)=0$ and a proper map $W:X^\vee\rightarrow \mathbb C$, where $W$ is called the superpotential.
\end{defn}

\begin{defn}\label{def-period}
Given a Landau-Ginzburg model and a choice of holomorphic $n$-form $\omega$, we define the \emph{periods} of the LG model \emph{relative to $W$ and $\omega$} to be the period functions associated to the varying fibers of the LG model obtained by integrating transcendental cycles across the $n$-form $\omega$.

\end{defn}
\begin{remark}
We expect that if $(X^\vee,W)$ is the mirror LG model of $X$, then the smooth fibers of $W$ should be mirror to generic anticanonical hypersurfaces in $X$. Since the fibers are Calabi-Yau manifolds, the choice of $\omega$ is well-defined up to multiplication by a holomorphic function on the base of the LG model. 
Scaling the holomorphic form by a function has the effect of scaling the periods by the same function.
\end{remark}

\begin{remark}
It is often the case that our LG models will depend on various deformation parameters. 
We clarify here that the relative period functions we consider are functions of both the base variable of the LG model and the deformation parameters. 
It will often be the case that we will treat this object as a \emph{deforming family of periods} in which case we will distinguish the base variable of the LG model from the deformation parameters. 
\end{remark}

\begin{remark}
As described in \cite{doran_calabi-yau_2017}, the geometry of an $M_n$-polarized family of $K3$ surfaces is determined by the associated generalized functional invariant. 
Thus, as long as one uses the pull-back of the chosen holomorphic $2$-form on $\mathcal{X}_n$ to compute the periods of the fibration, then the periods of the family will be precisely the pull-backs of the periods of $\mathcal{X}_n$ via the functional invariant. 
Later on in this paper, we will have to scale these period functions appropriately. 

More generally, the LG models that we  work with in this paper corresponding to higher-dimensional Calabi-Yau manifolds will be constructed via pull-back. 
Therefore, in these cases too, the periods of the LG model will be determined by pulling back by rational functions and scaling appropriately. 
\end{remark}

For Fano varieties, the classical period of the Minkowski polynomial associated to the LG model was introduced in \cite{CCGGK}. Note that the periods in \cite{CCGGK} are power series in one variable, which is essentially the base parameter of the LG model. We would like to point out that the periods in Definition \ref{def-period} can be specialized to the classical periods in \cite{CCGGK}.

\begin{example}
We consider the period for the LG model of $\mathbb P^3$ along with its smooth anticanonical $K3$ surface. The LG model can be written as the fiberwise compactification of the following.
The potential is
\begin{align*}
W:X^\vee\rightarrow \mathbb C\\
(x_1,x_2,x_3,y)\mapsto y,
\end{align*}
where $X^\vee$ is
\[
X^\vee=\left\{ (x_1,x_2,x_3,y)\in (\mathbb C^*)^4\left| x_1+x_2+x_3+\frac{q}{x_1x_2x_3}=y\right. \right\}.
\]
Performing a change of variables $x_i=y\tilde{x}_i$, we find that $X^\vee$ is birational to 
\[
\tilde{x}_1\tilde{x}_2\tilde{x}_3(\tilde{x}_1+\tilde{x}_2+\tilde{x}_3-1)+\frac{q}{y^4}=0.
\]
That is, $X^\vee$ is the pull-back of the family $\mathcal{X}_2$ via the generalized function invariant
\[
\lambda=\frac{q}{y^4}.
\]
We obtain periods for this LG model via pull-back. 
For example, the following expression is a period function with respect to the $2$-form obtained by pulling back the $2$-form on $\mathcal{X}_2$:

\[
f_0^{\mathbb P^3}(q,y):=f_{0}^{\mathcal{X}_2}(\lambda)=\sum_{d\geq 0}\frac{(4d)!}{(d!)^4}\left(\frac{q}{y^4}\right)^d.
\]
By setting $q=1$ and $1/y=t$, we recover the classical period
\[
f_0^{\mathbb P^3}(t)=\sum_{d\geq 0}\frac{(4d)!}{(d!)^4}t^{4d}.
\]
Although this is rather a trivial example, in general, the classical period can be obtained from the period defined in Definition \ref{def-period} in a similar way: set all the complex parameters to $1$ and set the base parameter $y$ of the LG model to $1/t$. 
\end{example}

\subsubsection{The iterative structure of periods}
In this section, we take a detour to explain how periods of a family of Calabi-Yau $n$-fold, fibered by Calabi-Yau $(n-1)$-fold, can also be computed using generalized functional invariants. Given a family of Calabi-Yau $n$-fold fibered by Calabi-Yau $(n-1)$-fold, the period of the family of Calabi-Yau $n$-fold is the residue integral (over a closed loop around $0$) of the pullback of the period of its internal fibration of the Calabi-Yau $(n-1)$-fold by the generalized functional invariant. This iterative structure has already studied in \cite{doran_calabi-yau_2015}.

The quintic mirror family $\mathcal{X}$ of Calabi-Yau threefolds is defined by
$$ \left\{x_1x_2x_3y(x_1+x_2+x_3+y-1)+\psi=0\right\}.$$
By setting
$$\tilde{x}_1=\frac{x_1}{1-y},\tilde{x}_2=\frac{x_2}{1-y},\ \tilde{x}_3=\frac{x_3}{1-y},\ \lambda=\frac{\psi}{y(1-y)^4},$$
the quintic mirror family is written in terms of the quartic mirror family. 
That is, for each $\psi\in\mathbb{P}^1-\{0,\frac{1}{5^5},\infty\}$, $\mathcal{X}_\psi$ is fibered by $M_2$-polarized $K3$ surfaces and the fibration is governed by the functional invariant
$$\lambda=\frac{\psi}{y(1-y)^4}.$$

As is shown in \cite{doran_calabi-yau_2015}*{Proposition 5.1}, we can calculate the periods of the quintic mirror family of Calabi-Yau threefold by integrating the ``relative'' periods corresponding to the fibration structure. 
We review some of the details below, but refer the readers to \cite{doran_calabi-yau_2015} for more. 

 First, note that:
\begin{align*}
    \frac{1}{y(1-y)}f_0^{\mathcal{X}_2}(\lambda)=&\frac{1}{y(1-y)}\sum_{d_1\geq 0}\frac{(4d_1)!}{(d_1!)^4}\left(\frac{\psi}{y(1-y)^4}\right)^{d_1}\\
    =&\sum_{d_1,d_2\geq 0}\frac{(4d_1)!}{(d_1!)^4}\left(\frac{\psi^{d_1}}{y^{d_1}}\right)\frac{(4d_1+d_2)!}{(4d_1)!d_2!}y^{d_2-1}\\
    =&\sum_{d_1,d_2\geq 0}\frac{(4d_1+d_2)!}{(d_1!)^4d_2!}\psi^{d_1}y^{d_2-d_1-1}.
\end{align*}
Here, we have used the fact that
$$\frac{1}{(1-y)^{k+1}}=\sum_{d=0}^\infty\frac{(d+k)!}{d!k!}y^d.$$
Next, if we integrate around a closed loop around $y=0$ and use the residue theorem, we find that
$$\frac{1}{2\pi i}\oint\frac{f_0^{\mathcal{X}_2}(\lambda)}{y(1-y)}=\sum_{d_1\geq 0}\psi^{d_1}\frac{(5d_1)!}{(d_1!)^5}.$$
The expression on the right-hand side is the well-known expression for the holomorphic period on the quintic mirror family of Calabi-Yau threefolds. 

We remark that the series expression above is only convergent for $\left|\frac{\psi}{y(1-y)^4}\right|<1$ and so the series formula present is only valid on the intersection of $\left|\frac{\psi}{y(1-y)^4}\right|<1$ and the cylinder $|y|<1$ (which is happily non-empty!).

In summary, we have constructed the standard holomorphic period of the quintic mirror family by integrating the $\psi$-dependent family of relative periods with respect to $y$. 
The scaling factor of $y(1-y)$ in front of the relative periods is present to ensure that the resulting integral corresponds precisely to the specific choice of holomorphic $3$-form that one normally makes when studying the quintic mirror.

\section{Tyurin degeneration of a conifold transition of quintic threefolds}\label{sec:coni-tran}

In this section, we consider the following Calabi-Yau threefold with two K\"ahler parameters: complete intersection of bidegrees (4,1) and (1,1) in $\mathbb P^4\times \mathbb P^1$. We denoted it by $\tilde {Q}_5$. It is a conifold transition of a quintic threefold in $\mathbb P^4$, see, for example, \cite{chialva_deforming_2008}. The Calabi-Yau threefold $\tilde {Q}_5$ admits a Tyurin degeneration
\[
\tilde Q_5\leadsto \tilde X_1\cup_{K3} \tilde X_2,
\]
where $\tilde X_1$ is a hypersurface of bidegree $(4,1)$ in $\mathbb P^3\times \mathbb P^1$ and $\tilde X_2$ is a complete intersection of bidegrees $(4,0), (1,1)$ in $\mathbb P^4\times \mathbb P^1$. Indeed, $\tilde X_1$ is the blow-up of $\mathbb P^3$ along the complete intersection of two quartic surfaces and $\tilde X_2$ is the blow-up of a quartic threefold $Q_4$ along the complete intersection of two hyperplanes (degree one hypersurfaces in $Q_4$).

\subsection{Generalized functional invariants}\label{sec:gen-fun-inv-conifold}
The mirrors of $\tilde Q_5$, $\tilde X_1$ and $\tilde X_2$ can be written down explicitly following \cite{Givental98}. The mirror $\tilde Q_5^\vee$ of $\tilde Q_5$ is the compactification of 
\[
\left\{(x_1,x_2,x_3,x_4,y)\in (\mathbb C^*)^5\left|
x_1+x_2+x_3+\frac{q_1}{x_1x_2x_3x_4}+\frac{q_0}{y}=1;
x_4+y =1\right.\right\}.
\]

The LG model for $\tilde X_1$ is 
\begin{align*}
W_1: \tilde X_1^\vee &\rightarrow \mathbb C\\
 (x_1,x_2,x_3,y_1)&\mapsto y_1,
\end{align*}
where $\tilde X_1^\vee$ is the fiberwise compactification of
\[
\left\{(x_1,x_2,x_3,y_1)\in(\mathbb C^*)^4\left|x_1+x_2+x_3+\frac{q_{1,1}}{x_1x_2x_3}+\frac{q_{0,1}}{y_1}=1\right.\right\}.
\]
The LG model for $\tilde X_2$ is
\begin{align*}
W_2: \tilde X_2^\vee&\rightarrow \mathbb C\\
 (x_1,x_2,x_3,y_2)&\mapsto  y_2,
\end{align*}
where $\tilde X_2^\vee$ is the fiberwise compactification of
\[
\left\{(x_1,x_2,x_3,y_2)\in(\mathbb C^*)^4\left|x_1+x_2+x_3+\frac{q_{1,2}}{x_1x_2x_3(1-q_{0,2}/y_2)}=1\right.\right\}.
\]

By performing an appropriate change of variables, we see that all three of these families are fibered by quartic mirror $K3$ surfaces.
This allows us to conclude that they are families of $M_2$-polarized $K3$ surfaces and read off their generalized functional invariants. 
For example, for the family $\tilde X_1^\vee$, we make the following change of variables:
\[
\tilde{x}_1=\frac{x_1}{1-q_{0,1}/y_1}, \tilde{x}_2=\frac{x_2}{1-q_{0,1}/y_1}, \tilde{x}_3=\frac{x_3}{1-q_{0,1}/y_1}, \lambda_1=\frac{q_{1,1}}{(1-q_{0,1}/y_1)^4}.
\]
This produces the quartic mirror family of $K3$ surfaces:
\[
\tilde{x}_1\tilde{x}_2\tilde{x}_3(\tilde{x}_1+\tilde{x}_2+\tilde{x}_3-1)+\lambda_1=0.
\]
We read off the generalized functional invariant for $\tilde X_1^\vee$ as $\lambda_1=\frac{q_{1,1}}{(1-q_{0,1}/y_1)^4}$. Similarly, we compute the generalized functional invariants for $\tilde{Q}_5^\vee$ and $\tilde X_2^\vee$ by making appropriate changes of variable to match with the quartic mirror family.
We obtain the following generalized functional invariants for $\tilde{Q}_5^\vee$, $\tilde X_1^\vee$ and $\tilde X_2^\vee$:

\[
\lambda=\frac{q_1}{(1-y)(1-q_0/y)^4}, \quad \lambda_1=\frac{q_{1,1}}{(1-q_{0,1}/y_1)^4}, \quad \lambda_2=\frac{q_{1,2}}{1-q_{0,2}/y_2}.
\]
Before proceeding, we map out the locations of the branch points and singular fibers. For $\tilde{Q}_5^\vee$ with functional invariant $\lambda=\frac{1}{(1-y)\left(1-\frac{q_0}{y}\right)^4}$, we have 
$$\lambda^{-1}(0)=\{0,\infty\},\ \lambda^{-1}(\infty)=\{1,q_0\},\ \lambda^{-1}\left(\frac{1}{256}\right)=\{\textrm{five points}\}.$$
Similarly, we find
$$\lambda_1^{-1}(0)=\{0\},\ \lambda_1^{-1}(\infty)=\{q_{0,1}\},\ \lambda_1^{-1}\left(\frac{1}{256}\right)=\{\textrm{four points}\}.$$
$$\lambda_2^{-1}(0)=\{0\},\ \lambda_2^{-1}(\infty)=\{q_{0,2}\},\ \lambda_2^{-1}\left(\frac{1}{256}\right)=\{\textrm{one point}\}.$$

The point $y_1=0$ is a cusp and corresponds to a semistable fiber of type $\textrm{III}$ with $34$ components; the fiber at $y_1=\infty$ is smooth; the four fibers located at the pre-image of $\frac{1}{256}$ are singular $K3$ surfaces containing a single $A_1$ singularity; the point $y_1=q_{0,1}$ has a type $\textrm{III}$ fiber. 
On the second LG model, the point $y_2=0$ is a cusp corresponding to a type $\textrm{III}$ fiber with $4$ components; the fiber $y_2=q_{0,2}$ has $31$ components; the fiber $y_2=\infty$ is smooth; the fiber over the pre-image of $\frac{1}{256}$ is a singular $K3$ surface with an $A_1$ singularity. 

The conifold transition itself has type $\textrm{III}$ fibers at $y=0$ and $y=\infty$ with $34$ and $4$ components respectively; the fiber $y=1$ has $31$ components; the fiber $y=q_0$ has a type $\textrm{III}$ fiber; the five fibers over the pre-images of $\frac{1}{256}$ are singular $K3$ surfaces with $A_1$ singularities. 

We glue the bases of the two LG models by making the following identification of variables:
\begin{align}\label{identification-coni-tran}
    q_1=q_{1,1}=q_{1,2},\quad q_0=q_{0,1},\quad y=y_1=q_{0,2}/y_2,
\end{align}
where the $q$-variables are naturally identified because they are coming from the same ambient space and the last equation is a natural identification based on the above analysis of the branch points and singular fibres.

\begin{proposition}
Under the identification (\ref{identification-coni-tran}), the following product relation holds among the functional invariants:
\begin{align}\label{equ:gen-fun-inv}
\frac{\lambda} {q_1}=\frac{\lambda_1}{q_{1,1}}\cdot\frac{\lambda_2}{q_{1,2}}.
\end{align}
\end{proposition}

\begin{remark}
The parameter $q_1$ above is a scaling parameter on the modular curve $X_0(2)^+$. 
Thus, one should think of equation \eqref{equ:gen-fun-inv} as saying that the generalized functional invariant of $\tilde{Q}_5^\vee$ is the product of the two functional invariants corresponding the LG models $\tilde{X}_1^\vee$ and $\tilde{X}_2^\vee$ after scaling appropriately. 
\end{remark}

This gluing formula for generalized functional invariants illustrates how the fibers of the two Landau-Ginzburg models are combined into the internal fibration on the Calabi-Yau threefold.  Before exploring the meaning of this gluing in terms of both fiberwise periods (of the $M_2$-polarized $K3$ surface fibers) and periods of the Calabi-Yau threefold itself, it is natural to ask whether the mirror to the conifold transition can be seen directly in terms of the generalized functional invariant $\lambda$.  The answer is a resounding "yes!".

\subsubsection{Mirror quintic as a limit}
Re-writing the functional invariant, we have
\[\lambda=\frac{yq_1}{(y-1)(y-q_0)^4}.\]
This corresponds to the family of $h^{2,1}=2$ Calabi-Yau threefolds that we are studying. 
In order to obtain the quintic mirror via a limit, we make a change of variables and take an appropriate limit.
First, make the transformation:
$$y=(1-q_0)\tilde{y}+q_0.$$
In this coordinate system, the generalized functional invariant transforms to
\[\lambda=\frac{q_1}{(1-q_0)^5}\frac{\left((1-q_0)\tilde{y}+q_0\right)}{\tilde{y}^4(\tilde{y}-1)}.\]
Now pull out a factor of $q_0$ from the numerator:
\[\lambda=\frac{q_1q_0}{(1-q_0)^5}\cdot\frac{\left(\frac{1-q_0}{q_0}\tilde{y}+1\right)}{\tilde{y}^4(\tilde{y}-1)}\]
Set $z_1=\frac{q_1}{(1-q_0)^4},\ z_2=z_1\frac{q_0}{1-q_0}=\frac{q_1q_0}{(1-q_0)^5}$. Then, the functional invariant above is re-written as 
\[\lambda=\frac{z_1\tilde{y}+z_2}{\tilde{y}(\tilde{y}-1)^4}.\]
In the limit $z_1\to0$, we obtain the quintic mirror with $z_2$ as the scaling parameter. 
In other words, we want $\frac{q_1}{(1-q_0)^4}$ to go to zero in such a way that $\frac{q_1q_0}{(1-q_0)^5}$ remains finite. 
This is exactly the limit that was considered in \cite{chialva_deforming_2008}.
The corresponding limits performed on the level of holomorphic periods are carried out in \cite{chialva_deforming_2008}.

\subsubsection{Topological Gluing}

\begin{figure}
    \centering
    \includegraphics[width=5cm]{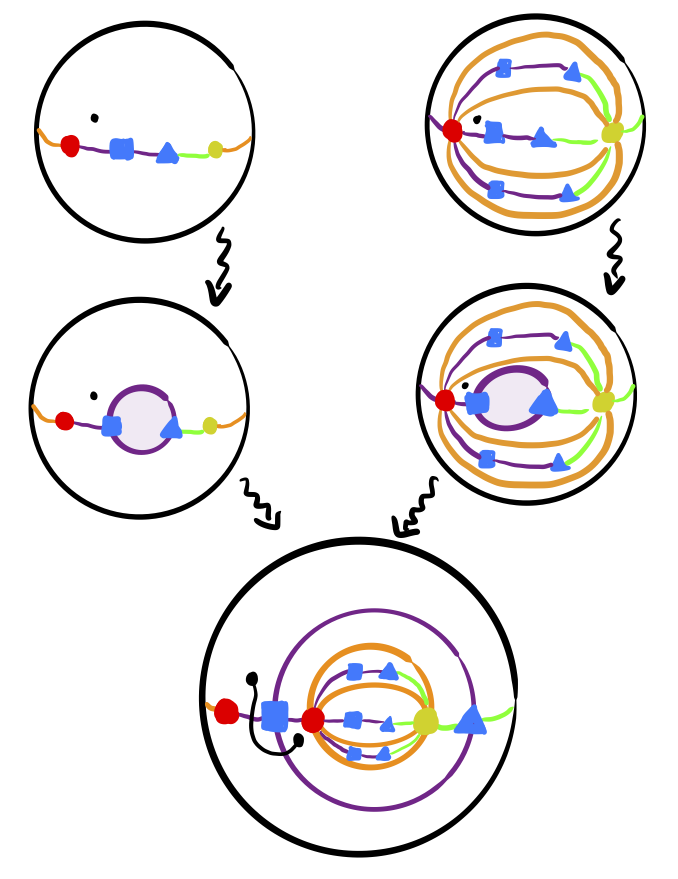}
    \caption{The topological gluing of functional invariants.}
    \label{fig:gluing}
\end{figure}

Before continuing, let us briefly explain how the multiplication of functional invariants described above relates to the topological gluing described in \cite{doran_hodge_2017}.
The maps $\lambda_1,\lambda_2,\lambda$ are all covers of the modular curve $X_0(2)^+$ branched over $0,\infty$, with $\lambda$ branching over two points $p_1,p_2$ depending on $z_1,z_2$. 
On the base curve $X_0(2)^+$, draw arcs joining $0$ to $p_1$, $p_1$ to $p_2$, $p_2$ to $\infty$ and $\infty$ to $0$. 
Then, the covers $\lambda_1,\lambda_2$ are determined by graphs drawn on the Riemann sphere. 
These graphs are depicted in Figure \ref{fig:gluing}, where the red dots are pre-images of $0$, the blue squares are pre-images of $p_1$, the blue triangles are pre-images of $p_2$, and the yellow dots are pre-images of $\infty$. 

The spheres at the top of Figure \ref{fig:gluing} correspond to $\lambda_1$ on the right and $\lambda_2$ on the left. 
Below, we have cut one of the edges joining a blue square to the adjacent blue triangle and opened it up. 
Finally, we glue the two spheres along these cuts. 
The resulting object is a sphere whose graph determines the cover corresponding to $\lambda$, which corresponds, in turn,  to the $h^{2,1}=2$ family of Calabi-Yau threefolds. 

In the limit as $z_1\to 0$, the blue squares degenerate into the red dots, which themselves coalesce to produce one red dot corresponding to the cusp of ramification index $5$. 
This cover corresponds to the quintic mirror. 

We can recognize the local system of the VHS corresponding to the $M_2$-polarized K3 surfaces as the gluing of the two individual VHSs corresponding to $\lambda_1,\lambda_2$. 
This follows from the fact that all three of these local systems are determined by the covers $\lambda,\lambda_1,\lambda_2$ respectively, but we can also be more explicit if we reference Figure \ref{fig:gluing}. 
Indeed, choose base points on each of the covering spheres, depicted as the black dots in Figure \ref{fig:gluing}, together with a basis of loops---note that the monodromy around the blue dots is trivial. 
Then, a basis of loops for the glued sphere is obtained by first taking the basis on the left sphere, and then ``dragging'' the basis of loops corresponding to the right sphere along the black path. 
Since there is no monodromy around the blue square, this is well-defined and the resulting global mondromy representation can be thought of as the concatenation of the original two monodromy representations. 
That is, if the monodromy tuple corresponding to the left sphere is $(g_1,g_2,g_3)$ and the monodromy truple corresponding to the right sphere is $(h_1,h_2,h_3,h_4,h_5,h_6)$, then the monodromy tuple of the gluing is $(g_1,g_2,g_3,h_1,h_2,h_3,h_4,h_5,h_6)$. 
Moreover, if we make these choices so that $g_3$ and $h_1$ correspond the monodromy around the red dots, then the limiting monodromy is obtained by simply multiplying $g_3$ and $h_1$ together.

\subsection{Relation among periods}
In Section \ref{sec:gen-fun-inv-conifold}, we saw how the generalized functional invariant for $\tilde{Q}_5$ arises by gluing two LG models and their corresponding generalized functional invariants. 
In this section, we explore what this means on the level of periods. 
\subsubsection{Periods for mirror quartic $K3$}
Consider the period associated to quartic mirror
\begin{align}
f_0^{K3}(q_1)=\sum_{d_1\geq 0}\frac{(4d_1)!}{(d_1!)^4}q_1^{d_1}.
\end{align}
Let $\theta_{q_1}=q_1\frac{\partial}{\partial q_1}$, then $f_0^{K3}(q_1)$ satisfies the following ODE
\begin{align}\label{ODE-K3}
(\theta_{q_1}^3-4q_1(4\theta_{q_1}+1)(4\theta_{q_1}+2)(4\theta_{q_1}+3))F(q_1)=0.    
\end{align}
The basis of solutions for Equation (\ref{ODE-K3}) can be taken as the coefficients of powers of $H$ in the following function: 
\begin{align}
    I^{K3}(q_1)=q_1^H\sum_{d_1\geq 0}^\infty q_1^{d_1}\frac{\prod_{k=1}^{4d_1}(4H+k)}{\prod_{k=1}^{d_1}(H+k)^4} \quad \mod H^3,
\end{align}
where $q_1^H=e^{H\log q_1}$.

\subsubsection{Periods for $\tilde Q_5^\vee$}

The holomorphic period for $\tilde Q_5$ can be found in \cite{chialva_deforming_2008}. It can also be obtained as a residue integral of the pullback of \[
\frac{1}{(1-y)(1-q_0/y)}f_0^{K3}(q_1)\]
by the generalized functional invariant
\[
\lambda=\frac{q_1}{(1-y)(1-q_0/y)^4}.
\]
We have
\begin{align*}
\frac{1}{(1-y)(1-q_0/y)}f_0^{K3}(\lambda)&=\frac{1}{(1-y)(1-q_0/y)}\sum_{d_1\geq 0}\frac{(4d_1)!}{(d_1!)^4}\left(\frac{q_1}{(1-y)(1-q_0/y)^4}\right)^{d_1}\\
&=\frac{1}{(1-y)}\sum_{d_1,d_{0,1}\geq 0}\frac{(4d_1)!}{(d_1!)^4}\left(\frac{q_1}{(1-q_0/y)^4}\right)^{d_1}\frac{(d_1+d_{0,1})!}{d_1!d_{0,1}!}(y)^{d_{0,1}}\\
&=\sum_{d_1,d_{0,1},d_{0,2}\geq 0}\frac{(4d_1)!}{(d_1!)^4}\left(q_1\right)^{d_1}\frac{(d_1+d_{0,1})!}{d_1!d_{0,1}!}\frac{(4d_1+d_{0,2})!}{(4d_1)!d_{0,2}!}(y)^{d_{0,1}}(q_0/y)^{d_{0,2}}.
\end{align*}

After taking a residue, we obtain the well-known formula for the holomorphic period of $\tilde{Q}_5^\vee$: 
\begin{align*}
f_0^{\tilde Q_5}(q_1,q_0)&=\frac{1}{2\pi i}\oint \frac{1}{(1-y)(1-q_0/y)}f_0^{K3}(\lambda) \frac{dy}{y}\\
&=\sum_{d_1,d_0\geq 0}\frac{(4d_1)!}{(d_1!)^4}\left(q_1\right)^{d_1}\frac{(d_1+d_{0})!}{d_1!d_{0}!}\frac{(4d_1+d_{0})!}{(4d_1)!d_{0}!}(q_0)^{d_{0}}\\
&=\sum_{d_1,d_0\geq 0}\frac{(4d_1+d_{0})!(d_1+d_{0})!}{(d_1!)^5(d_{0}!)^2}q_1^{d_1}q_0^{d_{0}}.
\end{align*}

\begin{remark}
These series will converge on the intersection $\left|\frac{q_1}{(1-y)(1-\frac{q_0}{y})}\right|<1$ and $|y|<1$. Similar convergence considerations will percolate in the work that follows.  We will avoid commenting further on domains of convergence unless there is likely to be ambiguity as to which domains we should be using. 
\end{remark}

The holomorphic period $f_0^{\tilde Q_5}(q_1,q_0)$ is the solution for the following system of PDEs.
\begin{align}
   &\left((\theta_{q_1})^5-q_1(4\theta_{q_1}+\theta_{q_0}+1)(4\theta_{q_1}+\theta_{q_0}+2)(4\theta_{q_1}+\theta_{q_0}+3)(4\theta_{q_1}+\theta_{q_0}+4)(\theta_{q_1}+\theta_{q_0}+1)\right)F(q_1,q_0)=0\\
  \notag &\left((\theta_{q_0})^2-q_0(4\theta_{q_1}+\theta_{q_0}+1)(\theta_{q_1}+\theta_{q_0}+1)\right)F(q_1,q_0)=0.
\end{align}

The basis of solutions for the system of PDEs is given by
\[
    I^{\tilde Q_5}(q_1,q_0)=q_1^H q_0^P\sum_{d_1,d_0\geq 0}^\infty q_1^{d_1}q_0^{d_0}\frac{\prod_{k=1}^{4d_1+d_0}(4H+P+k)\prod_{k=1}^{d_1+d_0}(H+P+k)}{\prod_{k=1}^{d_1}(H+k)^5\prod_{k=1}^{d_0}(P+k)^2},
\]
where $H$ and $P$ are hyperplane classes of $\mathbb P^4$ and $\mathbb P^1$ respectively; $I^{\tilde Q_5}(q_1,q_0)$ takes value in the ambient cohomology ring $\iota^*H^*(\mathbb P^4\times \mathbb P^1)$ with $\iota: \tilde Q_5\hookrightarrow \mathbb P^4\times \mathbb P^1$ being the inclusion.

\subsubsection{Periods for LG models}\label{sec:period-LG-example-1}

Following Section \ref{sec:period}, periods for $\tilde X_1^\vee$ can be constructed by pulling back periods of the mirror quartic via the generalized functional invariant map $\lambda_1$.
Consider the following scaled version of the pull-back of the holomorphic period:
\begin{align*}
f_0^{\tilde X_1}(q_{1,1},q_{0,1},y_1)=&\frac{1}{(1-q_{0,1}/y_1)}f_0^{K3}(\lambda_1)\\
=&\frac{1}{(1-q_{0,1}/y_1)}\sum_{d_1\geq 0}\frac{(4d_1)!}{(d_1!)^4}\left(\frac{q_{1,1}}{(1-q_{0,1}/y_1)^4}\right)^{d_1}\\
=&\sum_{d_1,d_0\geq 0}\frac{(4d_1)!(4d_1+d_0)!}{(d_1!)^4(4d_1)!d_0!}q_{1,1}^{d_1}(q_{0,1}/y_1)^{d_0}\\
=& \sum_{d_1,d_0\geq 0}\frac{(4d_1+d_0)!}{(d_1!)^4 d_0!}q_{1,1}^{d_1}(q_{0,1}/y_1)^{d_0}.
\end{align*}
Because $q_{0,1}$ and $y_1$ always appear at the same time with the form $q_{0,1}/y_1$, we will also denote by $f_0^{\tilde X_1}(q_{1,1},x)$ the holomorphic period for $\tilde X_1^\vee$ defined above:
\[
f_0^{\tilde X_1}(q_{1,1},x)=\sum_{d_1,d_0\geq 0}\frac{(4d_1+d_0)!}{(d_1!)^4 d_0!}q_{1,1}^{d_1}x^{d_0}.
\]
The holomorphic period $f_0^{\tilde X_1}(q_{1,1},x)$ is the holomorphic solution for the following system of PDEs.
\begin{align}
   &\left((\theta_{q_{1,1}})^4-q_{1,1}(4\theta_{q_{1,1}}+\theta_x+1)(4\theta_{q_{1,1}}+\theta_x+2)(4\theta_{q_{1,1}}+\theta_x+3)(4\theta_{q_{1,1}}+\theta_x+4)\right)F(q_{1,1},x)=0\\
  \notag &\left((\theta_x)-x(4\theta_{q_{1,1}}+\theta_x+1)\right)F(q_{1,1},x)=0,
\end{align}
where $F(q_{1,1},x)$ is a function of $q_{1,1}$ and $x$.

The basis of solutions for the system of PDEs is given by
\[
    I^{\tilde X_1}(q_{1,1},x)=q_{1,1}^H x^P\sum_{d_1,d_0\geq 0}^\infty q_{1,1}^{d_1}x^{d_0}\frac{\prod_{k=1}^{4d_1+d_0}(4H+P+k)}{\prod_{k=1}^{d_1}(H+k)^4\prod_{k=1}^{d_0}(P+k)}.
\]
We also write
\[
    I^{\tilde X_1}(q_{1,1},q_{0,1},y)=q_{1,1}^H (q_{0,1}/y)^P\sum_{d_1,d_0\geq 0}^\infty q_{1,1}^{d_1}(q_{0,1}/y)^{d_0}\frac{\prod_{k=1}^{4d_1+d_0}(4H+P+k)}{\prod_{k=1}^{d_1}(H+k)^4\prod_{k=1}^{d_0}(P+k)},
\]
where $H$ and $P$ are hyperplane classes of $\mathbb P^3$ and $\mathbb P^1$ respectively; $I^{\tilde X_1}$ takes value in the ambient cohomology ring $\iota^*H^*(\mathbb P^3\times \mathbb P^1)$ with $\iota: \tilde X_1\hookrightarrow \mathbb P^3\times \mathbb P^1$ being the inclusion.

Similarly, we can pull back the holomorphic period for the mirror quartic via the generalized functional invariant map $\lambda_2$ to obtain the holomorphic period for $\tilde X_2^\vee$:
\begin{align*}
\frac{1}{(1-q_{0,2}/y_2)}f_0^{K3}(\lambda_2)=&\frac{1}{(1-q_{0,2}/y_2)}\sum_{d_1,d_0\geq 0}\frac{(4d_1)!}{(d_1!)^4}\left(\frac{q_{1,2}}{1-q_{0,2}/y_2}\right)^{d_1}\\
=&\sum_{d_1,d_0\geq 0}\frac{(4d_1)!(d_1+d_0)!}{(d_1!)^4(d_1)!d_0!}q_{1,2}^{d_1}(q_{0,2}/y_2)^{d_0}\\
=& \sum_{d_1,d_0\geq 0}\frac{(4d_1)!(d_1+d_0)!}{(d_1!)^5 d_0!}q_{1,2}^{d_1}(q_{0,2}/y_2)^{d_0}.
\end{align*}

We define $f_0^{\tilde X_2}(q_{1,2},y)$, the holomorphic period for $\tilde X_2^\vee$, to be
\[
f_0^{\tilde X_2}(q_{1,2},y)=\sum_{d_1,d_0\geq 0}\frac{(4d_1)!(d_1+d_0)!}{(d_1!)^5 d_0!}q_{1,2}^{d_1}y^{d_0}.
\]
The holomorphic period $f_0^{\tilde X_2}(q_{1,2},y)$ is the solution for the following system of PDEs.
\begin{align}
   &\left((\theta_{q_{1,2}})^4-4q_{1,2}(4\theta_{q_{1,2}}+1)(4\theta_{q_{1,2}}+2)(4\theta_{q_{1,2}}+3)(\theta_{q_{1,2}}+\theta_y+1)\right)F(q_{1,2},y)=0\\
  \notag &\left((\theta_y)-y(\theta_{q_{1,2}}+\theta_y+1)\right)F(q_{1,2},y)=0,
\end{align}
where $F(q_{1,2},x)$ is a function of $q_{1,2}$ and $x$.

The basis of solutions for the system of PDEs is given by
\[
    I^{\tilde X_2}(q_{1,2},y)=q_{1,2}^H y^P\sum_{d_1,d_0\geq 0}^\infty q_{1,2}^{d_1}y^{d_0}\frac{\prod_{k=1}^{4d_1}(4H+k)\prod_{k=1}^{d_1+d_0}(H+P+k)}{\prod_{k=1}^{d_1}(H+k)^5\prod_{k=1}^{d_0}(P+k)},
\]
where $H$ and $P$ are hyperplane classes of $\mathbb P^4$ and $\mathbb P^1$ respectively; $I^{\tilde X_2}$ takes value in the ambient cohomology ring $\iota^*H^*(\mathbb P^4\times \mathbb P^1)$ with $\iota: \tilde X_2\hookrightarrow \mathbb P^4\times \mathbb P^1$ being the inclusion.

\begin{remark}
Let us address the scaling factors that appear in our period expressions. Ultimately, we will take the relative periods of $\tilde X_1^\vee$ and $\tilde X_2^\vee$ to construct the periods for $\tilde{Q}_5$. 
In order to make sure that the periods for $X_i^\vee$ are ``compatible'' with each other, we need to scale the periods appropriately. 
More precisely, the functions $f_0^{X_i}(q_{1,i},x)$ are $q_{1,i}$-dependent families of periods functions for the base variable $x$. 
The scaling factors are chosen to ensure that the characteristic exponents of the corresponding $q_{1,i}$-dependent families of ODEs at the singular points are the same for both $f_0^{\tilde X_1}(q_{1,1},x)$ and $f_0^{\tilde X_2}(q_{1,2},y)$. 
This should be thought of as a fine-tuning of the holomorphic form used to calculate the period integrals associated to each factor. 
\end{remark}

\subsubsection{Relations}

\begin{theorem}
The holomorphic periods of LG models can be glued together to form the holomorphic period of the mirror Calabi--Yau $\tilde Q_5^\vee$ with correction given by the holomorphic period of the mirror quartic $K3$ surface. More precisely, the relation is given by the Hadamard product
\[
f_0^{\tilde Q_5}(q_1,q_0)\star_{q_1} f_0^{K3}(q_1)=\frac{1}{2\pi i}\oint f_0^{\tilde X_1}(q_1,q_0,y)\star_{q_1} f_0^{\tilde X_2}(q_1,y)\frac{dy}{y},
\]
where $\star_{q_1}$ means the Hadamard product with respect to the variable $q_1$.
\end{theorem}

\begin{proof}
This is a straightforward computation: 
\begin{align*}
    &\frac{1}{2\pi i}\oint f_0^{\tilde X_1}(q_1,q_0,y)\star_{q_1} f_0^{\tilde X_2}(q_1,y)\frac{dy}{y}\\
    =&\frac{1}{2\pi i}\oint \sum_{d_1,d_{0,1},d_{0,2}\geq 0}\frac{(4d_1+d_{0,1})!}{(d_1!)^4 d_{0,1}!} \frac{(4d_1)!(d_1+d_{0,2})!}{(d_1!)^5 d_{0,2}!}q_1^{d_1}(q_0/y)^{d_{0,1}}(y)^{d_{0,2}}\frac{dy}{y}\\
    =& \sum_{d_1,d_0\geq 0}\frac{(4d_1+d_0)!}{(d_1!)^4 d_0!} \frac{(4d_1)!(d_1+d_0)!}{(d_1!)^5 d_0!}q_1^{d_1}q_0^{d_0}\\
    =& f_0^{\tilde Q_5}(q_1,q_0)\star_{q_1} f_0^{K3}(q_1).
\end{align*}
\end{proof}

\begin{remark}
Recall that, we have the following identification of variables:
\[
q_1=q_{1,1}=q_{1,2}, \quad q_0=q_{0,1}, \quad y=y_1=q_{0,2}/y_2.
\]
Gluing the LG models on the level of periods means taking the residue integral for the periods of LG models over $y$, where $y$ is the base parameter of the internal fibration of $\tilde Q_5^\vee$. The period of $\tilde Q_5^\vee$ are also computed by residue integral (over the parameter $y$) of period of the Calabi-Yau family in one dimensional lower via the generalized functional invariants. 
\end{remark}

One can also write down a Hadamard product relation among the Picard-Fuchs operators. 
Moreover, we have the following identity among the bases of solutions to the Picard-Fuchs equations.
\begin{theorem}
We have the following Hadamard product relation
\begin{align}\label{Hadamard-I-function-conifold}
    I^{\tilde Q_5}(q_1,q_0)\star_{q_1} I^{K3}(q_1)=\frac{1}{2\pi i}\oint I^{\tilde X_1}(q_1,q_0,y)\star_{q_1} I^{\tilde X_2}(q_1,y)\frac{dy}{y}.
\end{align}
In the Hadamard product $*_{q_1}$, we treat $\log q_1$ as a variable that is independent from $q_1$. Alternative, one may consider $\bar{I}(q_1,q_0):=I(q_1,q_0)/q_1^Hq_0^P$ and write Hadamard product relation for $\bar{I}$.
\end{theorem}

More explicitly, one can write an identity for each coefficient of $H^aP^b$ for $0\leq a \leq 2$, $0\leq b \leq 1$. Let $f_{i,j}^{\tilde Q_5}(q_1,q_0)$ be the coefficient of $H^iP^j$ for $I^{\tilde Q_5}(q_1,q_0)$. We still write $f_0^{\tilde Q_5}(q_1,q_0)$ for the coefficient of $H^0P^0$, since it is simply the holomorphic period. Similarly, for the coefficients of $I^{K3}(q_1)$, $I^{\tilde X_1}(q_1,q_0,y)$ and $I^{\tilde X_2}(q_1,y)$. We have the following identity for the coefficient of $H$:
\begin{align*}
&f_{1,0}^{\tilde Q_5}(q_1,q_0)\star_{q_1} f_{0}^{K3}(q_1)+f_{0}^{\tilde Q_5}(q_1,q_0)\star_{q_1} f_1^{K3}(q_1)\\
=&\frac{1}{2\pi i}\oint \left(f_{1,0}^{\tilde X_1}(q_1,q_0,y)\star_{q_1} f_0^{\tilde X_2}(q_1,y)+f_{0}^{\tilde X_1}(q_1,q_0,y)\star_{q_1} f_{1,0}^{\tilde X_2}(q_1,y)\right)\frac{dy}{y}.
\end{align*}
For the coefficient of $H^2$, we have
\begin{align*}\label{identity-H-2}
&f_{2,0}^{\tilde Q_5}(q_1,q_0)\star_{q_1} f_{0}^{K3}(q_1)+f_{1,0}^{\tilde Q_5}(q_1,q_0)\star_{q_1} f_{1}^{K3}(q_1)+f_{0}^{\tilde Q_5}(q_1,q_0)\star_{q_1} f_2^{K3}(q_1)\\
\notag =&\frac{1}{2\pi i}\oint \left(f_{2,0}^{\tilde X_1}(q_1,q_0,y)\star_{q_1} f_0^{\tilde X_2}(q_1,y)+f_{1,0}^{\tilde X_1}(q_1,q_0,y)\star_{q_1} f_{1,0}^{\tilde X_2}(q_1,y)+f_{0}^{\tilde X_1}(q_1,q_0,y)\star_{q_1} f_{2,0}^{\tilde X_2}(q_1,y)\right)\frac{dy}{y}.
\end{align*}

\begin{remark}
We want to point out that $I^{\tilde X_1}$ and $I^{\tilde X_2}$ are not exactly the $I$-functions for the corresponding relative Gromov-Witten invariants in Section \ref{sec:mirror-theorem-pairs}. We will use $I^{(\tilde X_1,K3)}$ and $I^{(\tilde X_2,K3)}$ to denote the relative $I$-functions which are more complicated than $I^{\tilde X_1}$ and $I^{\tilde X_2}$. Nevertheless, the information of $I^{\tilde X_1}$ and $I^{\tilde X_2}$ can be extracted from the relative $I$-functions.
\end{remark}

\section{Tyurin degeneration of quintic threefolds}\label{sec:quintic}

\subsection{Blow-up along $\mathbb P^3$}
We consider the Tyurin degeneration of a quintic threefold $Q_5$ into a quartic threefold $Q_4$ and the blow-up $\on{Bl}_C\mathbb P^3$ of $\mathbb P^3$ along complete intersection center of degrees $4$ and $5$ hypersurfaces, that is,
\[
Q_5\leadsto Q_4\cup_{K3}\on{Bl}_C\mathbb P^3,
\]
where $C$ is the complete intersection of degrees $4$ and $5$ hypersurfaces in $\mathbb P^3$, and $K3$ is the common anticanonical $K3$ hypersurface. We would like to write down the LG models for $Q_4$ and $\on{Bl}_C\mathbb P^3$ and glue them to the mirror family of $Q_5$. Since $\on{Bl}_C\mathbb P^3$ can be written as a hypersurface in a toric variety, we can write down the LG models following Givental \cite{Givental98}. For the rest of this section, we write $X_1:=Q_4$ and $\tilde X_2:=\on{Bl}_C\mathbb P^3$.

The LG model for $X_1=Q_4$ is
\begin{align*}
W_1:X_1^\vee\rightarrow \mathbb C\\
(x_1,x_2,x_3,y_1)\mapsto y_1,
\end{align*}
where $X_1^\vee$ is the fiberwise compactification of
\[
\left\{(x_1,x_2,x_3,y_1)\in (\mathbb C^*)\left| x_1+x_2+x_3+\frac{q_{1,1}}{x_1x_2x_3y_1}=1\right.\right\}.
\]

Following \cite{CCGK}*{Section E}, the quasi-Fano variety $\tilde X_2=\on{Bl}_C(\mathbb P^3)$ can be constructed as a hypersurface of degree (4,1) in the toric variety $\mathbb P(\mathcal O_{\mathbb P^3}(-1)\oplus \mathcal O_{\mathbb P^3})$.

The LG model of $\tilde X_2=\on{Bl}_C(\mathbb P^3)$ is defined as follows
\begin{align*}
W_2:\tilde X_2^\vee\rightarrow \mathbb C\\
(x_1,x_2,x_3,y_2)\rightarrow y_2,
\end{align*}
where $\tilde X_2^\vee$ is the fiberwise compactification of
\[
\left\{(x_1,x_2,x_3,y_2)\in (\mathbb C^*)^4\left|x_1+x_2+x_3+\frac{q_{0,2}}{y_2}+\frac{q_{1,2}y_2}{x_1x_2x_3}=1\right.\right\}.
\]
Finally, the mirror quintic family $Q_5^\vee$ is defined by the following equation
\[
x_1x_2x_3y(x_1+x_2+x_3+y-1)+q_1=0.
\]

Similar to the computation in Section \ref{sec:coni-tran}, we obtain the generalized functional invariants for $Q_5^\vee$, $X_1^\vee$ and $\tilde X_2^\vee$ respectively: 
\begin{align}
    \lambda=\frac{q_1}{y(1-y)^4}, \quad \lambda_1=\frac{q_{1,1}}{y_1}, \quad \lambda_2=\frac{q_{1,2}y_2}{(1-q_{0,2}/y_2)^4}.
\end{align}
We set
\begin{align}\label{identification-Q_5-1}
q_1=q_{1,1}=q_{1,2}y_2, \quad y=y_1=q_{0,2}/y_2.
\end{align}
The matching of singular fibers works similarly to the previous section with singular fibers on one LG model being glued to smooth fibers of the other. 
We have
\begin{proposition}
Under the identification (\ref{identification-Q_5-1}), the following product relation holds among the functional invariants:
\begin{align}\label{equ:gen_fun_inv2}
\frac{\lambda}{q_1}=\frac{\lambda_1}{q_{1,1}}\cdot\frac{\lambda_2}{q_{1,2}}.
\end{align}
\end{proposition}

\begin{remark}
Similar to Equation \eqref{equ:gen-fun-inv}, Equation \eqref{equ:gen_fun_inv2} says that the the functional invariant for $Q_5^\vee$ is equal to the product of the functional invariants of the LG models $X_1^\vee$ and $\tilde{X}_2^\vee$ after scaling.
\end{remark}

The period for $Q_5^\vee$ is
\[
f_0^{Q_5}(q_1)=\sum_{d_1\geq 0}q_1^{d_1}\frac{(5d_1)!}{(d_1!)^5}.
\]
The period for $X_1^\vee$ is
\[
f_0^{X_1}(q_{1,1},y_1)=f_0^{K3}(\lambda_1)=\sum_{d_1\geq 0}\frac{(4d_1)!}{(d_1!)^4}\left(\frac{q_{1,1}}{y_1}\right)^{d_1}.
\]
The period for $\tilde X_2^\vee$ is
\begin{align*}
f_0^{\tilde X_2}(q_{1,2}y_2,q_{0,2}/y_2)=&\frac{1}{1-y}f_0^{K3}(\lambda_2)\\
=&\frac{1}{1-q_{0,2}/y_2}\sum_{d_1\geq 0}\frac{(4d_1)!}{(d_1!)^4}\left(\frac{q_{1,2}y_2}{(1-q_{0,2}/y_2)^4}\right)^{d_1}\\
=& \sum_{d_1,d_0\geq 0}\frac{(4d_1)!}{(d_1!)^4}(q_{1,2}y_2)^{d_1}\frac{(4d_1+d_0)!}{(4d_1)!d_0!}(q_{0,2}/y_2)^{d_0}\\
=& \sum_{d_1,d_0\geq 0}\frac{(4d_1+d_0)!}{(d_1)!^4d_0!}(q_{1,2}y_2)^{d_1}(q_{0,2}/y_2)^{d_0}.
\end{align*}

We may rewrite the period for $\tilde X_2^\vee$ as
\[
f_0^{\tilde X_2}(q_1,y)=\sum_{d_1,d_0\geq 0}\frac{(4d_1+d_0)!}{(d_1)!^4d_0!}(q_1)^{d_1}(y)^{d_0}.
\]

Note that $q_1=q_{1,1}=q_{1,2}y_2$ and $y=y_1=q_{0,2}/y_2$ when we glue the LG models.
\begin{theorem}
We have the following Hadamard product relation:
\[
f_0^{Q_5}(q_1)\star_{q_1} f_0^{K3}(q_1)=\frac{1}{2\pi i}\oint f_0^{X_1}(q_1,y)\star_{q_1} f_0^{\tilde X_2}(q_1,y)\frac{dy}{y}.
\]
\end{theorem}
\begin{proof}
This is a straightforward computation.
\end{proof}
The Picard-Fuchs operators and the bases of solutions to Picard-Fuchs equations are related in a similar way. Let
\[
I^{Q_5}(q_1)=e^{H\on{log}q_1}\sum_{d\geq 0}\left(\frac{\prod_{k=1}^{5d} (5H+k)}{\prod_{k=1}^d(H+k)^{5}}\right)q_1^d;
\]
\begin{align*}
    I^{X_1}(q_1,y)=(q_1/y)^H\sum_{d_1\geq 0}^\infty (q_1/y)^{d_1}\frac{\prod_{k=1}^{4d_1}(4H+k)}{\prod_{k=1}^{d_1}(H+k)^4};
\end{align*}
\[
    I^{\tilde X_2}(q_1,y)=q_1^H y^P\sum_{d_1,d_0\geq 0}^\infty q_1^{d_1}y^{d_0}\frac{\prod_{k=1}^{4d_1+d_0}(4H+P+k)}{\prod_{k=1}^{d_1}(H+k)^4\prod_{k=1}^{d_0}(P+k)}.
\]
\begin{theorem}
We have the following Hadamard product relation
\begin{align}
    I^{Q_5}(q_1)\star_{q_1} I^{K3}(q_1)=\frac{1}{2\pi i}\oint I^{X_1}(q_1,y)\star_{q_1} I^{\tilde X_2}(q_1,y)\frac{dy}{y},
\end{align}
where we set $P=H$. In the Hadamard product $*_{q_1}$, we treat $\log q_1$ as a variable that is independent from $q_1$. Alternative, we can consider $\bar{I}(q_1):=I(q_1)/q_1^H$ and write Hadamard product relation for $\bar{I}$.
\end{theorem}
\subsection{Blow-up along the quartic threefold.}

We consider the Tyurin degeneration of $Q_5$ into $\mathbb P^3$ and the blow-up of $Q_4$ along the complete intersection $C$ of hypersurfaces of degrees $1$ and $5$:
\[
Q_5\leadsto \mathbb P^3\cup_{K3}\on{Bl}_C(Q_4).
\]
For the rest of this section we write $X_1:=\mathbb P^3$ and $\tilde X_2:=\on{Bl}_C(Q_4)$.

The LG model for $X_1=\mathbb P^3$ is given by the fiberwise compactification of
\begin{align*}
W: (\mathbb C^*)^3&\rightarrow \mathbb C\\
 (x_1,x_2,x_3)&\mapsto x_1+x_2+x_3+\frac {q_{1,1}}{x_1x_2x_3}.
\end{align*}
It can be rewritten as follows.
The potential is
\begin{align*}
W_1:X_1^\vee\rightarrow \mathbb C\\
(x_1,x_2,x_3,y_1)\mapsto y_1,
\end{align*}
where $X_1^\vee$ is the fiberwise compactification of
\[
\left\{ (x_1,x_2,x_3,y_1)\in (\mathbb C^*)^4\left| x_1+x_2+x_3+\frac{q_{1,1}}{x_1x_2x_3}=y_1\right. \right\}.
\]

The blown-up variety $\tilde X_2=\on{Bl}_C(Q_4)$ can be realized as a complete intersection of degrees $(1,1)$ and $(4,0)$ in the toric variety $\mathbb P(\mathcal O_{\mathbb P^4}(-4)\oplus \mathcal O_{\mathbb{P}^4})\rightarrow \mathbb P^4$. The LG model of $\tilde X_2=\on{Bl}_C(Q_4)$ can be written as follows:

The potential is
\begin{align*}
W_2:\tilde X_2^\vee\rightarrow \mathbb C\\
(x_1,x_2,x_3,y_2)\mapsto y_2,
\end{align*}
where $\tilde X_2^\vee$ is the fiberwise compactification of
\[
\left\{ (x_1,x_2,x_3,y_2)\in (\mathbb C^*)^4\left| x_1+x_2+x_3+\frac{y_2^4q_{1,2}}{x_1x_2x_3(1-q_{0,2}/y_2)}=1\right. \right\}.
\]

The generalized functional invariants for $Q_5^\vee$, $X_1^\vee$ and $\tilde X_2^\vee$ are
\[
    \lambda=\frac{q_1}{y(1-y)^4}, \quad \lambda_1=\frac{q_{1,1}}{y_1^4}, \quad \lambda_2=\frac{q_{1,2}y_2^4}{(1-q_{0,2}/y_2)}.
\]
For $\lambda$, we consider the change of variable $y\mapsto 1-y$, then we have
\[
\lambda=\frac{q_1}{y^4(1-y)}, \quad \lambda_1=\frac{q_{1,1}}{y_1^4}, \quad \lambda_2=\frac{q_{1,2}y_2^4}{(1-q_{0,2}/y_2)}.
\]

We consider the following identification among variables
\begin{align}\label{identification-Q_5-2}
q_1=q_{1,1}=q_{1,2}y_2^4, \quad y=y_1=q_{0,2}/y_2.
\end{align}

\begin{proposition}
Under the identification (\ref{identification-Q_5-2}), we have the relation among generalized functional invariants
\[
\frac{\lambda}{q_1}=\frac{\lambda_1}{q_{1,1}}\cdot \frac{\lambda_2}{q_{1,2}}.
\]
\end{proposition}

The holomorphic period for $X_1^\vee$ is
\[
f_0^{X_1}(q_{1,1},y_1)=f_0^{K3}(\lambda_1)=\sum_{d_1\geq 0}\frac{(4d_1)!}{(d_1!)^4}\left(\frac{q_{1,1}}{y_1^4}\right)^{d_1}.
\]
The holomorphic period for $\tilde X_2^\vee$ is
\begin{align*}
f_0^{\tilde X_2}(q_{1,2}y_2^4,q_{0,2}/y_2)=\frac{1}{1-y}f_0^{K3}(\lambda_2)=&\frac{1}{1-q_{0,2}/y_2}\sum_{d_1\geq 0}\frac{(4d_1)!}{(d_1!)^4}\left(\frac{q_{1,2}y_2^4}{(1-q_{0,2}/y_2)}\right)^{d_1}\\
=& \sum_{d_1,d_0\geq 0}\frac{(4d_1)!}{(d_1!)^4}\frac{(d_1+d_0)!}{d_1!d_0!}(q_{1,2}y_2^4)^{d_1}(q_{0,2}/y_2)^{d_0}.
\end{align*}
We can rewrite the period for $\tilde X_2^\vee$ as
\[
f_0^{\tilde X_2}(q_1,y)=\sum_{d_1,d_0\geq 0}\frac{(4d_1)!}{(d_1!)^4}\frac{(d_1+d_0)!}{d_1!d_0!}(q_1)^{d_1}(y)^{d_0}.
\]

We have the Hadamard product relation among periods.
\begin{theorem} The following relation holds for holomorphic periods:
\[
f_0^{Q_5}(q_1)\star_{q_1} f_0^{K3}(q_1)=\frac{1}{2\pi i}\oint f_0^{X_1}(q_1,y)\star_{q_1} f_0^{\tilde X_2}(q_1,y)\frac{dy}{y}.
\]
\end{theorem}

The Picard-Fuchs operators and the bases of solutions to Picard-Fuchs equations are related in a similar way. Let
\begin{align*}
    I^{X_1}(q_1,y)=(q_1/y^4)^H\sum_{d_1\geq 0}^\infty (q_1/y^4)^{d_1}\frac{\prod_{k=1}^{4d_1}(4H+k)}{\prod_{k=1}^{d_1}(H+k)^4};
\end{align*}
\[
    I^{\tilde X_2}(q_1,y)=q_1^H y^P\sum_{d_1,d_0\geq 0}^\infty q_1^{d_1}y^{d_0}\frac{\prod_{k=1}^{4d_1}(4H+k)\prod_{k=1}^{d_1+d_0}(H+P+k)}{\prod_{k=1}^{d_1}(H+k)^5\prod_{k=1}^{d_0}(P+k)}.
\]
\begin{theorem}
We have the following Hadamard product relation
\begin{align}
    I^{Q_5}(q_1)\star_{q_1} I^{K3}(q_1)=\frac{1}{2\pi i}\oint I^{X_1}(q_1,y)\star_{q_1} I^{\tilde X_2}(q_1,y)\frac{dy}{y},
\end{align}
where we set $P=4H$. In the Hadamard product $\star_{q_1}$, we treat $\log q_1$ as a variable that is independent from $q_1$. Alternative, we can consider $\bar{I}(q_1):=I(q_1)/q_1^H$ and write Hadamard product relation for $\bar{I}$.
\end{theorem}
\section{Tyurin degeneration of Calabi-Yau complete intersections in toric varieties}

\subsection{Set-up}\label{sec:set-up}
Let $X$ be a Calabi-Yau complete intersection in a toric variety $Y$ defined by a generic section of $E=L_0\oplus L_1\oplus\cdots \oplus  L_s$, where each $L_l$ is a nef line bundle. Let $\rho_l=c_1(L_l)$, then $-K_Y=\sum_{l=0}^s c_1(L_l)$. Let $s_l\in H^0(Y,L_l)$ be generic sections determining $X$. A refinement of the nef partition with respect to $L_0$ is given by two nef line bundles $L_{0,1}, L_{0,2}$ such that $L_0=L_{0,1}\otimes L_{0,2}$. Let $\rho_{0,1}=c_1(L_{0,1})$ and $\rho_{0,2}=c_1(L_{0,2})$. We have two quasi-Fano varieties $X_1$ and $X_2$ defined by sections of $E_1=L_{0,1}\oplus L_1\oplus\cdots \oplus  L_s$ and $E_2=L_{0,2}\oplus L_1\oplus\cdots \oplus  L_s$ respectively. Let $s_{0,1}\in H^0(Y,L_{0,1})$ and $s_{0,2}\in H^0(Y,L_{0,2})$ be the generic sections determining $X_1$ and $X_2$. For our construction of a Tyurin degeneration, we blow up $X_2$ along $X_2\cap \{s_0=0\}\cap \{s_{0,1}=0\}$. We write the blown-up variety as $\tilde{X}_2$. Two quasi-Fano varieties $X_1$ and $\tilde X_2$ intersect along $X_0$ which is a Calabi-Yau complete intersection in the toric variety $Y$ defined by a generic section of $L_{0,1}\oplus L_{0,2}\oplus L_1\oplus\cdots \oplus  L_s$.

We can again realize $\tilde{X}_2$ as a complete intersection in a toric variety following \cite{CCGK}*{Section E}. Indeed, it is a hypersurface in the total space of $\pi:\mathbb P_{X_2}(\mathcal O\oplus i^*L_{0,2}^{-1})\rightarrow X_2$ defined by a generic section of the line bundle $\pi^*i^*L_{0,1}\otimes \mathcal O(1)$, where $i:X_2\hookrightarrow Y$ is the inclusion map. In other words, it is a complete intersection in the toric variety $\mathbb P_Y(\mathcal O\oplus L_{0,2}^{-1})$ given by a generic section of $\pi^*L_{0,2}\oplus \pi^*L_1\oplus\cdots \oplus  \pi^*L_s \oplus (\pi^*L_{0,1}\otimes \mathcal O(1))$ where we use the same $\pi$ for the projection of $\mathbb P_Y(\mathcal O\oplus L_{0,2}^{-1})$ to the base $Y$. Then $-K_{\tilde{X}_2}=(P-\pi^*i^*\rho_{0,2})|_{\tilde X_2}$ by the adjunction formula, where $P=c_1(\mathcal O_{P_Y(\mathcal O\oplus L_{0,2}^{-1})}(1))$.

\begin{rmk}
In dimension one, there is no codimension-two subvarieties to blow up. However, we may use the same geometric construction for blow-ups as above to construct a one dimensional subvariety in another toric variety. Then the results in dimension one are parallel to the results in higher dimensions.
\end{rmk}

Let $p_1,\ldots,,p_r\in H^2(Y,\mathbb Z)$ be a nef integral basis. We write the toric divisors as
\[
D_j=\sum_{i=1}^r m_{ij}p_i,\quad 1\leq j\leq m, 
\]
for some $m_{ij}$.

Let $X$ be the Calabi-Yau complete intersection in the toric variety $Y$. The nef partition of the toric divisors gives a partition of the variables $x_1,\ldots, x_m$ into $s+1$ groups. Let $F_l(x)$ be the sum of $x_i$ in each group $l=0,\ldots,s$. Following \cite{Givental98}, the Hori-Vafa mirror $X^\vee$ of $X$ is Calabi-Yau compactification of 
\[
\left\{(x_1,\ldots,x_m)\in (\mathbb C^*)^m\left| \prod_{j=1}^m x_j^{m_{ij}}=q_i, i=1,\ldots,r; F_l(x)=1,l=0,\ldots,s\right.\right\}.
\]
Note that $q_i$, $1\leq i\leq r$, are complex parameters for $X^\vee$ from the ambient toric variety $Y$.

\subsection{Blowing-up both quasi-Fanos}\label{sec:blow-up-2}

We first consider the case when we blow-up both quasi-Fanos. Then, the Calabi-Yau and quasi-Fanos, denoted by $\tilde{X},\tilde{X}_1$ and $\tilde{X}_2$ respectively, become complete intersections in the toric variety $Y\times \mathbb P^1$. Let $\pi_1$ and $\pi_2$ be the projection of $Y\times \mathbb P^1$ onto $Y$ and $\mathbb P^1$ respectively. Then 
\begin{itemize}
    \item  $\tilde X$ is the complete intersection in $Y\times \mathbb P^1$ defined by generic sections of $\pi_1^*L_1$, $\ldots$, $\pi_1^*L_s$, $\pi^*_1 L_{0,1}\otimes \pi_2^*\mathcal O_{\mathbb P^1}(1)$ and $\pi^*_1 L_{0,2}\otimes \pi_2^*\mathcal O_{\mathbb P^1}(1)$; 
    \item $\tilde X_1$ is the complete intersection in $Y\times \mathbb P^1$ defined by generic sections of $\pi_1^*L_1$, $\ldots$, $\pi_1^*L_s$, $\pi^*_1 L_{0,1}$ and $\pi^*_1 L_{0,2}\otimes \mathcal \pi_2^*\mathcal O_{\mathbb P^1}(1)$; 
    \item $\tilde X_2$ is the complete intersection in $Y\times \mathbb P^1$ defined by generic sections of $\pi_1^*L_1$, $\ldots$, $\pi_1^*L_s$, $\pi^*_1 L_{0,1}\otimes \mathcal \pi_2^*\mathcal O_{\mathbb P^1}(1)$ and $\pi^*_1 L_{0,2}$;
    \item $X_0$ is the complete intersection in $Y$ defined by generic sections of $L_{0,1}$, $L_{0,2}$, $L_1$, $\cdots$, $L_s$.
\end{itemize}

The mirror for $\tilde{X}$ is the compactification of
\begin{align*}
&\left\{(x_1,\ldots,x_m,y_1,y_2)\in (\mathbb C^*)^{m+2}\left| \prod_{j=1}^m x_j^{m_{ij}}=q_i, i=1,\ldots,r; y_1y_2=q_0;\right.\right.\\ & \left.\vphantom{\prod_{j=1}^m} F_l(x)=1,l=1,\ldots,s;F_{0,1}(x)+y_1=1,F_{0,2}(x)+y_2=1\right\}\\
=&\left\{(x_1,\ldots,x_m,y)\in (\mathbb C^*)^{m+1}\left| \prod_{j=1}^m x_j^{m_{ij}}=q_i, i=1,\ldots,r;\right.\right.\\ &\left. \vphantom{\prod_{j=1}^m} F_l(x)=1,l=1,\ldots,s;F_{0,1}(x)+y=1,F_{0,2}(x)+\frac{q_0}{y}=1\right\},
\end{align*}
where $F_{0,1}(x)$ and $F_{0,2}(x)$ correspond to the refinement of the nef partition with respect to $L_0$. In particular, we have $F_{0,1}(x)+F_{0,2}(x)=F_0(x)$.

The mirror for $\tilde{X}_1$ is a LG model
\begin{align*}
W:\tilde{X}_1^\vee &\rightarrow \mathbb C^*\\
(x_1,\ldots,x_m,y_1)&\mapsto y_1,
\end{align*}
where $\tilde{X}_1^\vee$ is the fiberwise compactification of
\begin{align*}
&\left\{(x_1,\ldots,x_m,y_1)\in (\mathbb C^*)^{m+1}\left| \prod_{j=1}^m x_j^{m_{ij}}=q_{i,1}, i=1,\ldots,r;\right.\right.\\ &\left.\vphantom{\prod_{j=1}^m} F_l(x)=1,l=1,\ldots,s;F_{0,1}(x)=1,F_{0,2}(x)+\frac{q_{0,1}}{y_1}=1\right\}.
\end{align*}

Similarly, the mirror for $\tilde{X}_2$ is a LG model
\begin{align*}
W:\tilde{X}_2^\vee &\rightarrow \mathbb C^*\\
(x_1,\ldots,x_m,y_2)&\mapsto y_2,
\end{align*}
where $\tilde{X}_2^\vee$ is the fiberwise compactification of
\begin{align*}
&\left\{(x_1,\ldots,x_m,y_2)\in (\mathbb C^*)^{m+1}\left| \prod_{j=1}^m x_j^{m_{ij}}=q_{i,2}, i=1,\ldots,r;\right.\right.\\ &\left.\vphantom{\prod_{j=1}^m} F_l(x)=1,l=1,\ldots,s;F_{0,1}(x)+\frac{q_{0,2}}{y_2}=1,F_{0,2}(x)=1\right\}.
\end{align*}

The $\tilde X_1$ and $\tilde X_2$ intersect along $X_0$ which is a Calabi-Yau variety in one dimensional lower. The mirror for $X_0$ is $X_0^\vee$ which is the compactification of 
\begin{align*}
\left\{(x_1,\ldots,x_m)\in (\mathbb C^*)^m\left| \prod_{j=1}^m x_j^{m_{ij}}=q_i, i=1,\ldots,r;\right.\right. \\
\left.\vphantom{\prod_{j=1}^m} F_l(x)=1,l=1,\ldots,s; F_{0,1}(x)=1, F_{0,2}(x)=1\right\}.
\end{align*}

\subsubsection{Generalized functional invariants}\label{sec:gen-fun-inv}

We can compute the generalized functional invariants from the mirrors. Let $F_{0,1}(x)=x_{j_1}+\ldots,x_{j_a}$ and $F_{0,2}(x)=x_{j_{a+1}}+\ldots,x_{j_{a+b}}$. For $\tilde{X}^\vee$, we have $F_{0,1}(x)+y=1$ and $F_{0,2}(x)+\frac{q_0}{y}=1$. For $y\neq 1$ and $\frac{q_0}{y} \neq 1$, the following change of variables can give us $X_0^\vee$.
\[
\bar{x}_{j_k}=\frac{x_{j_k}}{1-y}, k=1,\ldots,a;\quad \bar{x}_{j_k}=\frac{x_{j_k}}{1-q_0/y}, k=a+1,\ldots,a+b; 
\]
and 
\begin{align}\label{gen-fun-inv-tilde-X}
\lambda_i=\frac{q_i}{\prod_{k=1}^a(1-y)^{m_{ij_k}}\prod_{k=a+1}^{a+b}(1-q_0/y)^{m_{ij_k}}}, i=1\ldots, r.
\end{align}
Then, $\lambda=(\lambda_1,\ldots,\lambda_r)$ is called the generalized functional invariant for $\tilde{X}^\vee$.

The generalized functional invariants for $\tilde{X}_1^\vee$ and $\tilde{X}_2^\vee$ can be computed in a similar way. Indeed, we have
\begin{align}\label{gen-fun-inv-tilde-X-1}
\lambda_{i,1}=\frac{q_{i,1}}{\prod_{k=a+1}^{a+b}(1-q_{0,1}/y_1)^{m_{ij_k}}}, i=1\ldots, r;
\end{align}
\begin{align}\label{gen-fun-inv-tilde-X-2}
\lambda_{i,2}=\frac{q_{i,2}}{\prod_{k=1}^a(1-q_{0,2}/y_2)^{m_{ij_k}}}, i=1\ldots, r.
\end{align}

We have the following identification of the variables:
\begin{align}\label{identification-tilde-X}
y=y_1=q_{0,2}/y_2,\quad q_0=q_{0,1}, \quad q_i=q_{i,1}=q_{i,2}, 1\leq i\leq r.
    \end{align}

\begin{proposition}
Under the identification (\ref{identification-tilde-X}), the generalized functional invariants satisfy the product relation
\[
\frac{\lambda_i}{q_i}=\frac{\lambda_{i,1}}{q_{i}}\frac{\lambda_{i,2}}{q_{i}}, \quad i=1,\ldots, r.
\]
\end{proposition}

\begin{remark}\label{remark-specialization}
To specialize to Section \ref{sec:coni-tran}, we set $Y=\mathbb P^4$. Then $m=5$, $r=1$, $l=0$, and $m_{1j}=1$ for $j=1,\ldots, 5$. We also have $F_{0,1}=x_4$ and $F_{0,2}=x_1+x_2+x_3+x_5$, where $x_5=\frac{Q}{x_1x_2x_3x_4}$. Then the generalized functional invariants specialize to the generalized functional invariants in Section \ref{sec:coni-tran}.
\end{remark}

\subsubsection{Periods}

Let $q=(q_1,\ldots,q_r)$ and $\lambda=(\lambda_1,\ldots,\lambda_r)$. The holomorphic period for $X_0^\vee$ is
\[
f_0^{X_0}(q)=\sum_{\substack{d\in H_2(Y;\mathbb Z)\\ \forall j, \langle D_j, d\rangle\geq 0}}\frac{\prod_{l=1}^s \langle \rho_l, d\rangle !}{\prod_{j=1}^m \langle D_j, d\rangle!}\langle \rho_{0,1}, d\rangle!\langle \rho_{0,2}, d\rangle ! q^d,
\]
where $q^d=q_1^{\langle p_1, d\rangle}\cdots q_r^{\langle p_r, d\rangle}$.
We can compute the holomorphic period for $\tilde{X}^\vee$ via generalized functional invariants:
\begin{align*}
&\frac{1}{(1-y)(1-q_0/y)}f_0^{X_0}(\lambda)\\
=&\sum_{\substack{d\in H_2(Y;\mathbb Z)\\ \forall j, \langle D_j, d\rangle \geq 0}}\frac{\prod_{l=1}^s \langle \rho_l, d\rangle!}{\prod_{j=1}^m \langle D_j, d\rangle !}\langle \rho_{0,1}, d\rangle!\langle\rho_{0,2}, d\rangle ! \frac{q^d}{(1-y)^{1+\langle \rho_{0,1}, d\rangle }(1-q_0/y)^{1+\langle \rho_{0,2}, d\rangle}}\\
=&\sum_{\substack{d\in H_2(Y;\mathbb Z)\\ \forall j, \langle D_j, d\rangle \geq 0}}\sum_{d_{0,1},d_{0,2}\geq 0}\frac{\prod_{l=1}^s \langle \rho_l, d\rangle !}{\prod_{j=1}^m \langle D_j, d\rangle !}\frac{(d_{0,1}+\langle \rho_{0,1}, d\rangle )!(d_{0,2}+\langle \rho_{0,2}, d\rangle )!q^d}{d_{0,1}!d_{0,2}!}y^{d_{0,1}}(q_0/y)^{d_{0,2}}.
\end{align*}
Then, we obtain the holomorphic period for $\tilde{X}^\vee$:
\begin{align*}
f_0^{\tilde{X}}(q,q_0)&=\frac{1}{2\pi i}\oint \frac{1}{(1-y)(1-q_0/y)}f_0^{X_0}(\lambda)\frac{dy}{y}\\
&=\sum_{\substack{d\in H_2(Y;\mathbb Z)\\ \forall j, \langle D_j, d\rangle \geq 0}}\sum_{d_{0}\geq 0}\frac{\prod_{l=1}^s \langle \rho_l, d\rangle !}{\prod_{j=1}^m \langle D_j, d\rangle !}\frac{(d_{0}+\langle \rho_{0,1}, d\rangle) !(d_{0}+\langle \rho_{0,2}, d\rangle)!q^d}{(d_{0}!)^2}q_0^{d_0}.
\end{align*}
The holomorphic periods for $\tilde{X}_1^\vee$ and $\tilde{X}_2^\vee$ are defined as the pullback of the corresponding generalized functional invariants, we have
\[
f_0^{\tilde{X}_1}(q,q_{0,1},y_1)=\sum_{\substack{d\in H_2(Y;\mathbb Z)\\ \forall j, \langle D_j, d\rangle \geq 0}}\sum_{d_{0}\geq 0}\frac{\prod_{l=1}^s \langle \rho_l, d\rangle!}{\prod_{j=1}^m \langle D_j, d\rangle !}\frac{\langle \rho_{0,1}, d\rangle !(d_{0}+\langle \rho_{0,2}, d\rangle )!q^d}{d_{0}!}(q_{0,1}/y_1)^{d_0};
\]
\[
f_0^{\tilde{X}_2}(q,y)=\sum_{\substack{d\in H_2(Y;\mathbb Z)\\ \forall j, \langle D_j, d\rangle \geq 0}}\sum_{d_{0}\geq 0}\frac{\prod_{l=1}^s \langle \rho_l, d\rangle !}{\prod_{j=1}^m \langle D_j, d\rangle !}\frac{(d_{0}+\langle \rho_{0,1}, d\rangle )!\langle \rho_{0,2}, d\rangle !q^d}{d_{0}!}(y)^{d_0}.
\]

Then, by direct computation, we have 
\begin{theorem}
The following Hadamard product relation holds
\[
f_0^{\tilde{X}}(q,q_0)\star_q f_0^{X_0}(q)=\frac{1}{2\pi i}\oint f_0^{\tilde{X}_1}(q,q_0,y)\star_q f_0^{\tilde{X}_2}(q,y)\frac{dy}{y}.
\]
\end{theorem}
Let $\on{NE}(Y)\subset H_2(Y,\mathbb R)$ be the cone generated by effective curves and $\on{NE}(Y)_{\mathbb Z}:=\on{NE}(Y)\cap (H_2(Y,\mathbb Z)/\on{tors})$. We have the Hadamard product relation among Picard-Fuchs equations. 
\begin{theorem}
The Hadamard product relation among the bases of solutions to Picard-Fuchs equations is
\[
I^{\tilde{X}}(q,q_0)\star_q I^{X_0}(q)=\frac{1}{2\pi i}\oint I^{\tilde{X}_1}(q,q_0/y)\star_q I^{\tilde{X}_2}(q,y)\frac{dy}{y},
\]
where
\begin{align*}
I^{\tilde X}(q,q_0)=&e^{\sum_{i=0}^r p_i\log q_i}\sum_{d\in \on{NE}(Y)_{\mathbb Z},d_0\geq 0}\prod_{j=1}^m \left(\frac{\prod_{k=-\infty}^0(D_j+k)}{\prod_{k=-\infty}^{\langle D_j,d\rangle}(D_j+k)}\right)\left(\prod_{l=1}^s \prod_{k=1}^{\langle \rho_l, d\rangle}(\rho_l+k)\right)\\
&\cdot\frac{\left(\prod_{k=1}^{\langle \rho_{0,1}, d\rangle+d_0}(\rho_{0,1}+p_0+k)\prod_{k=1}^{\langle \rho_{0,2},d\rangle+d_0}(\rho_{0,2}+p_0+k)\right)}{\prod_{k=1}^{d_0}(p_0+k)^2}q_0^{d_0}q^d;
\end{align*}
\begin{align*}
I^{X_0}(q)=&e^{\sum_{i=1}^r p_i\log q_i}\sum_{d\in \on{NE}(Y)_{\mathbb Z}}\prod_{j=1}^m \left(\frac{\prod_{k=-\infty}^0(D_j+k)}{\prod_{k=-\infty}^{\langle D_j,d\rangle}(D_j+k)}\right)\left(\prod_{l=1}^s \prod_{k=1}^{\langle \rho_l, d\rangle}(\rho_l+k)\right)\\
&\cdot\left(\prod_{k=1}^{\langle \rho_{0,1}, d\rangle }(\rho_{0,1}+k)\prod_{k=1}^{\langle \rho_{0,2}, d\rangle}(\rho_{0,2}+k)\right)q^{d};
\end{align*}
\begin{align*}
I^{\tilde X_1}(q,q_0)=e^{\sum_{i=0}^r p_i\log q_i}\sum_{d\in \on{NE}(Y)_{\mathbb Z},d_0\geq 0}&\prod_{j=1}^m \left(\frac{\prod_{k=-\infty}^0(D_j+k)}{\prod_{k=-\infty}^{\langle D_j,d\rangle}(D_j+k)}\right)\left(\prod_{l=1}^s \prod_{k=1}^{\langle \rho_l,d\rangle}(\rho_l+k)\right)\\
&\cdot\frac{\left(\prod_{k=1}^{\langle \rho_{0,1}, d\rangle }(\rho_{0,1}+k)\prod_{k=1}^{\langle \rho_{0,2}, d\rangle +d_0}(\rho_{0,2}+p_0+k)\right)}{\prod_{k=1}^{d_0}(p_0+k)}q_0^{d_0}q^d;
\end{align*}
\begin{align*}
I^{\tilde X_2}(q,q_0)=e^{\sum_{i=0}^r p_i\log q_i}\sum_{d\in \on{NE}(Y)_{\mathbb Z},d_0\geq 0}&\prod_{j=1}^m \left(\frac{\prod_{k=-\infty}^0(D_j+k)}{\prod_{k=-\infty}^{\langle D_j,d\rangle}(D_j+k)}\right)\left(\prod_{l=1}^s \prod_{k=1}^{\langle \rho_l, d\rangle}(\rho_l+k)\right)\\
&\cdot\frac{\left(\prod_{k=1}^{\langle \rho_{0,1}, d\rangle+d_0}(\rho_{0,1}+p_0+k)\prod_{k=1}^{\langle \rho_{0,2}, d\rangle}(\rho_{0,2}+k)\right)}{\prod_{k=1}^{d_0}(p_0+k)}q_0^{d_0}q^d.
\end{align*}
In the Hadamard product $*_q$, we treat $\log q_i$ as a variable that is independent from $q_i$. Alternative, we can consider $\bar{I}(q):=I(q)/(\sum p_i \log q_i)$ and write Hadamard product relation for $\bar{I}$.
\end{theorem}

\begin{remark}
The case when we blow-up both quasi-Fanos can be considered as a special case of the Tyurin degeneration when blowing-up occurs on only one of the quasi-Fanos. Indeed, we can directly consider the Tyurin degeneration for $\tilde{X}$ in $Y\times \mathbb P^1$ given by a refinement of the nef partition $(\pi_1^*L_{0,2}\otimes \pi_2^*\mathcal O_{\mathbb P^1}(1)) \oplus L_1\oplus\cdots \oplus  L_s \oplus (\pi_1^*L_{0,1}\otimes \pi_2^*\mathcal O_{\mathbb P^1}(1))$ with respect to either $\pi_1^*L_{0,2}\otimes \pi_2^*\mathcal O_{\mathbb P^1}(1))$ or $\pi_1^*L_{0,1}\otimes \pi_2^*\mathcal O_{\mathbb P^1}(1))$.  
\end{remark}

\subsection{General case}\label{sec:general-case}

We consider the case when we blow-up one of the quasi-Fanos. We assume that we blow-up $X_2$. The blown-up variety is denoted by $\tilde{X}_2$. The mirrors for $X$ and $X_0$ are described in Section \ref{sec:set-up}. The mirror for $X_1$ is the LG model
\begin{align*}
W:X_1^\vee &\rightarrow \mathbb C^*\\
(x_1,\ldots,x_m)&\mapsto F_{0,2}(x),
\end{align*}
where $X_1^\vee$ is the fiberwise compactification of
\begin{align*}
&\left\{(x_1,\ldots,x_m)\in (\mathbb C^*)^{m}\left| \prod_{j=1}^m x_j^{m_{ij}}=q_{i,1}, i=1,\ldots,r;\right.\right.\\ &\left.\vphantom{\prod_{j=1}^m} F_l(x)=1,l=1,\ldots,s;F_{0,1}(x)=1\right\}.
\end{align*}

$\tilde{X}_2$ is a complete intersection in the toric variety $\mathbb P_Y(\mathcal O\oplus L_{0,2}^{-1})$ given by a generic section of $L_{0,2}\oplus L_1\oplus\cdots \oplus  L_s \oplus (\pi^*L_{0,1}\otimes \mathcal O(1))$. The LG model for $\tilde{X}_2$ is
\begin{align*}
W:\tilde{X}_2^\vee &\rightarrow \mathbb C^*\\
(x_1,\ldots,x_m,y_2)&\mapsto y_2,
\end{align*}
where $\tilde{X}_2^\vee$ is the fiberwise compactification of
\begin{align*}
&\left\{(x_1,\ldots,x_m,y_2)\in (\mathbb C^*)^{m+1}\left| \left(\prod_{j=1}^m x_j^{m_{ij}}\right)y_2^{-\sum_{k=a+1}^{a+b}m_{ij_k}}=q_{i,2}, i=1,\ldots,r;\right.\right.\\ &\left.\vphantom{\prod_{j=1}^m} F_l(x)=1,l=1,\ldots,s;F_{0,1}(x)+\frac{q_{0,2}}{y_2}=1,F_{0,2}(x)=1 \right\}.
\end{align*}

\subsubsection{Generalized functional invariants}
Then, we can compute the generalized functional invariants. Recall that $X^\vee$ is the compactification of
\[
\left\{(x_1,\ldots,x_m)\in (\mathbb C^*)^m\left| \prod_{j=1}^m x_j^{m_{ij}}=q_i, i=1,\ldots,r; F_l(x)=1,l=0,\ldots,s\right.\right\}.
\]
Set $F_{0,2}(x)=y$, then $F_{0,1}(x)=1-y$. Following the computation in Section \ref{sec:gen-fun-inv}, the generalized functional invariant for $X^\vee$ is
\begin{align}\label{gen-fun-inv-X}
\lambda_i=\frac{q_i}{\prod_{k=1}^a(1-y)^{m_{ij_k}}\prod_{k=a+1}^{a+b}(y)^{m_{ij_k}}}, i=1\ldots, r.
\end{align}
For $X_1^\vee$, we also set $F_{0,2}(x)=y_1$, then the generalized functional invariant is
\begin{align}\label{gen-fun-inv-X-1}
\lambda_{i,1}=\frac{q_{i,1}}{\prod_{k=a+1}^{a+b}y_1^{m_{ij_k}}}, i=1\ldots, r.
\end{align}
The generalized functional invariant for $\tilde{X}_2$ is
\begin{align}\label{lambda-i-2}
\lambda_{i,2}=\frac{q_{i,2}}{y_2^{-\sum_{k=a+1}^{a+b}m_{ij_k}}\prod_{k=1}^{a}(1-q_{0,2}/y_2)^{m_{ij_k}}}, i=1\ldots, r.
\end{align}

We set
\begin{align}\label{identification-X}
q_i=q_{i,1}=q_{i,2}y_2^{\sum_{k=a+1}^{a+b}m_{ij_k}}, i=1,\ldots,r, \quad y=y_1=q_{0,2}/y_2.
\end{align}
\begin{proposition}
Under the identification (\ref{identification-X}), the generalized functional invariants satisfy the product relation
\[
\frac{\lambda_i}{q_i}=\frac{\lambda_{i,1}}{q_{i}}\frac{\lambda_{i,2}}{q_{i}}, \quad i=1,\ldots, r.
\]
\end{proposition}

\begin{remark}
One can specialize to Section \ref{sec:quintic} following Remark \ref{remark-specialization}.
\end{remark}

\begin{remark}
Section \ref{sec:blow-up-2} can be viewed as a special case of the current section as follows. The ambient toric variety in Section \ref{sec:blow-up-2} is $Y\times \mathbb P^1$. Let $p_1,\ldots,,p_r\in H^2(Y,\mathbb Z)$ be a nef integral basis and $p_0\in H^2(\mathbb P^1,\mathbb Z)$ is the hyperplane class in $\mathbb P^1$. We use the same notation $p_i$ to denote the pullbacks of $p_i$ to $H^2(Y\times \mathbb P^1,\mathbb Z)$. We can write the toric divisors as
\[
D_j=\sum_{i=1}^r m_{ij}p_i+0p_0,\quad 1\leq j\leq m,
\]
for some $m_{ij}$. And
\[
D_{m+1}=D_{m+2}=\sum_{i=1}^r 0p_i+ p_0.
\]
In other words, 
\[
D_j=\sum_{i=0}^r m_{ij}p_i,\quad 1\leq j\leq m+2,
\]
where $m_{ij}=0$ for $1 \leq j \leq m$ and $i=r+1$; $m_{ij}=0$ for $j=m+1, m+2$ and $1 \leq i \leq m$.
The Tyurin degeneration is given by the refinement of the nef partition with respect to the part $\pi_1^*L_{0,1}\otimes \pi_2^*\mathcal O_{\mathbb P^1}(1)$. The refinement produces two new parts $F_{0,1}(x)=x_{j_1}+\ldots,x_{j_a}$ and $F_{0,2}=y$. The generalized functional invariants for $X^\vee$ and $X_1^\vee$ are already computed in Section \ref{sec:blow-up-2}. The generalized functional invariant for $\tilde X_2^\vee$ is Equation (\ref{lambda-i-2}). In this special case, the factor $y_2^{-\sum_{k=a+1}^{a+b}m_{ij_k}}$ in the denominator of (\ref{lambda-i-2}) is actually $y_2^0=1$, because there is only one element in $F_{0,2}$ for $i\in\{1,\ldots,r\}$ and $j\in\{m+1,m+2\}$. So the generalized functional invariant for $\tilde X_2^\vee$ specialize to the one in Section \ref{sec:blow-up-2}. Hence, we recover the generalized functional invariants in Section \ref{sec:blow-up-2}. 

Note that in Section \ref{sec:blow-up-2}, we consider $X_0$ as a complete intersection in $Y$, while, in the current section, we consider it as a complete intersection in $Y\times \mathbb P^1$. This results in a slightly different expressions of its mirror $X_0^\vee$. Therefore, generalized functional invariants for $X^\vee$ and $X_1^\vee$ are changed accordingly. Nevertheless, the product relation among generalized functional invariants always holds.
\end{remark}

\begin{remark}
Note that 
\[
\rho_{0,2}=c_1(L_{0,2})=\sum_{k=a+1}^{a+b} D_{j_k}=\sum_{k=a+1}^{a+b}\sum_{i=1}^r m_{ij_k}p_i=\sum_{i=1}^r\left(\sum_{k=a+1}^{a+b} m_{ij_k}\right)p_i.
\]
The first Chern class of the normal bundle of the anticanonical divisor $X_0$ in $X_1$ is $i^*\rho_{0,2}\in H^2(X_0)$, where $i:X_0\hookrightarrow Y$ is the inclusion. It gives an explanation for the factor $y_2^{\sum_{k=a+1}^{a+b}m_{ij_k}}$ in Equation (\ref{identification-X}). This factor does not appear in Section \ref{sec:blow-up-2} because the normal bundle is trivial.
\end{remark}

\subsubsection{Periods}

The holomorphic period for $X_0^\vee$ is
\[
f_0^{X_0}(q)=\sum_{\substack{d\in H_2(Y;\mathbb Z)\\ \forall j, \langle D_j, d\rangle\geq 0}}\frac{\prod_{l=1}^s \langle \rho_l, d\rangle !}{\prod_{j=1}^m \langle D_j, d\rangle !}\langle \rho_{0,1}, d \rangle !\langle \rho_{0,2}, d\rangle ! q^d.
\]
We can compute the holomorphic period for $X^\vee$ via generalized functional invariants.
\begin{align*}
&\frac{1}{(1-y)}f_0^{X_0}(\lambda)\\
=&\sum_{\substack{d\in H_2(Y;\mathbb Z)\\ \forall j, \langle D_j, d\rangle \geq 0}}\frac{\prod_{l=1}^s \langle \rho_l, d\rangle !}{\prod_{j=1}^m \langle D_j, d\rangle !}\langle \rho_{0,1}, d\rangle !\langle \rho_{0,2}, d\rangle ! \frac{q^d}{(1-y)^{1+\langle \rho_{0,1}, d\rangle }y^{\langle \rho_{0,2}, d\rangle}}\\
=&\sum_{\substack{d\in H_2(Y;\mathbb Z)\\ \forall j, \langle D_j, d\rangle \geq 0}}\sum_{d_{0,1}\geq 0}\frac{\prod_{l=1}^s \langle \rho_l, d\rangle !}{\prod_{j=1}^m \langle D_j, d\rangle !}\langle \rho_{0,1}, d\rangle !\langle \rho_{0,2}, d\rangle !\frac{(d_{0,1}+\langle \rho_{0,1}, d\rangle) !q^d}{d_{0,1}!\langle \rho_{0,1}, d\rangle !y^{\langle \rho_{0,2}, d\rangle }}y^{d_{0,1}}.
\end{align*}
We obtain the holomorphic period for $X^\vee$:
\begin{align*}
f_0^{X}(q)&=\frac{1}{2\pi i}\oint \frac{1}{(1-y)}f_0^{X_0}(\lambda)\frac{dy}{y}\\
&=\sum_{\substack{d\in H_2(Y;\mathbb Z)\\ \forall j, \langle D_j, d\rangle \geq 0}}\frac{\prod_{l=1}^s \langle \rho_l, d\rangle !}{\prod_{j=1}^m \langle D_j, d\rangle !}(\langle \rho_{0,1}, d\rangle+\langle \rho_{0,2}, d\rangle)!q^d.\\
&=\sum_{\substack{d\in H_2(Y;\mathbb Z)\\ \forall j, \langle D_j, d\rangle \geq 0}}\frac{\prod_{l=0}^s \langle\rho_l, d\rangle !}{\prod_{j=1}^m \langle D_j, d\rangle !}q^d.
\end{align*}
The holomorphic periods for $X_1^\vee$ and $\tilde{X}_2^\vee$ are computed in a similar way, we have
\[
f_0^{X_1}(q,y_1)=\sum_{\substack{d\in H_2(Y;\mathbb Z)\\ \forall j, \langle D_j, d\rangle \geq 0}}\frac{\prod_{l=1}^s \langle \rho_l, d\rangle !}{\prod_{j=1}^m \langle D_j, d\rangle !}\langle \rho_{0,1}, d\rangle !\langle \rho_{0,2}, d\rangle !\frac{q^d}{(y_1)^{\langle \rho_{0,2}, d\rangle}};
\]
\[
f_0^{\tilde{X}_2}(q,y)=\sum_{\substack{d\in H_2(Y;\mathbb Z)\\ \forall j, \langle D_j, d\rangle \geq 0}}\sum_{d_{0}\geq 0}\frac{\prod_{l=1}^s \langle \rho_l, d\rangle !}{\prod_{j=1}^m \langle D_j, d\rangle !}\frac{(d_{0}+\langle \rho_{0,1}, d\rangle) !\langle \rho_{0,2}, d\rangle !q^d}{d_{0}!}(y)^{d_0}.
\]
Recall that $y=y_1=q_{0,2}/y_2$. Then we have 
\begin{theorem}\label{thm:gluing-period-general}
The following Hadamard product relation holds
\[
f_0^{X}(q)\star_q f_0^{X_0}(q)=\frac{1}{2\pi i}\oint f_0^{X_1}(q,y)\star_q f_0^{\tilde{X}_2}(q,y)\frac{dy}{y}.
\]
\end{theorem}
Similarly, we have the Hadamard product relation among Picard-Fuchs operators. Let
\begin{align*}
I^{X}(q)=e^{\sum_{i=0}^r p_i\log q_i}\sum_{d\in \on{NE}(Y)_{\mathbb Z}}&\prod_{j=1}^m \left(\frac{\prod_{k=-\infty}^0(D_j+k)}{\prod_{k=-\infty}^{\langle D_j,d\rangle}(D_j+k)}\right)\left(\prod_{l=0}^s \prod_{k=1}^{\langle \rho_l, d\rangle }(\rho_l+k)\right)q^d;
\end{align*}
\begin{align*}
I^{X_0}(q)=e^{\sum_{i=1}^r p_i\log q_i}\sum_{d\in \on{NE}(Y)_{\mathbb Z}}&\prod_{j=1}^m \left(\frac{\prod_{k=-\infty}^0(D_j+k)}{\prod_{k=-\infty}^{\langle D_j,d\rangle}(D_j+k)}\right)\left(\prod_{l=1}^s \prod_{k=1}^{\langle \rho_l, d\rangle }(\rho_l+k)\right)\\
&\cdot\left(\prod_{k=1}^{\langle \rho_{0,1}, d\rangle }(\rho_{0,1}+k)\prod_{k=1}^{\langle \rho_{0,2}, d\rangle }(\rho_{0,2}+k)\right)q^{d};
\end{align*}
\begin{align*}
I^{X_1}(q,y)=e^{\sum_{i=1}^r p_i\log \lambda_{i,1}}\sum_{d\in \on{NE}(Y)_{\mathbb Z},d_0\geq 0}&\prod_{j=1}^m \left(\frac{\prod_{k=-\infty}^0(D_j+k)}{\prod_{k=-\infty}^{\langle D_j,d\rangle}(D_j+k)}\right)\left(\prod_{l=1}^s \prod_{k=1}^{\langle \rho_l, d\rangle }(\rho_l+k)\right)\\
&\cdot\left(\prod_{k=1}^{\langle \rho_{0,1}, d\rangle }(\rho_{0,1}+k)\prod_{k=1}^{\langle \rho_{0,2}, d\rangle }(\rho_{0,2}+k)\right)\frac{q^d}{y^{\langle \rho_{0,2}, d\rangle }},
\end{align*}
where $\lambda_{i,1}=\frac{q_{i,1}}{\prod_{k=a+1}^{a+b}y_1^{m_{ij_k}}}, i=1\ldots, r$;

\begin{align*}
I^{\tilde X_2}(q,y)=e^{\sum_{i=1}^r p_i\log q_i+p_0\log y}\sum_{d\in \on{NE}(Y)_{\mathbb Z},d_0\geq 0}&\prod_{j=1}^m \left(\frac{\prod_{k=-\infty}^0(D_j+k)}{\prod_{k=-\infty}^{\langle D_j,d\rangle}(D_j+k)}\right)\left(\prod_{l=1}^s \prod_{k=1}^{\langle \rho_l, d\rangle }(\rho_l+k)\right)\\
&\cdot\frac{\left(\prod_{k=1}^{\langle \rho_{0,1}, d\rangle +d_0}(\rho_{0,1}+p_0+k)\prod_{k=1}^{\langle \rho_{0,2}, d\rangle }(\rho_{0,2}+k)\right)}{\prod_{k=1}^{d_0}(p_0+k)}y^{d_0}q^d.
\end{align*}
\begin{theorem}\label{thm:general-case-I-function}
The Hadamard product relation among the bases of solutions to Picard-Fuchs equations is
\[
I^{X}(q)\star_q I^{X_0}(q)=\frac{1}{2\pi i}\oint I^{X_1}(q,y)\star_q I^{\tilde{X}_2}(q,y)\frac{dy}{y},
\]
where we set $p_0=\rho_{0,2}$. In the Hadamard product $*_q$, we treat $\log q_i$ as a variable that is independent from $q_i$. Alternative, we can consider $\bar{I}(q):=I(q)/(\sum p_i \log q_i)$ and write Hadamard product relation for $\bar{I}$.
\end{theorem}

\section{Gromov-Witten invariants}
In this section, we discuss how the relation among periods is related to the A-model data: Gromov-Witten invariants. The Tyurin degeneration naturally relates absolute Gromov-Witten invariants of a Calabi-Yau variety and the relative Gromov-Witten invariants of the quasi-Fano varieties via the degeneration formula \cite{Li01}, \cite{Li02} and \cite{LR}. The periods for the mirror Calabi-Yau families are related to absolute Gromov-Witten invariants of Calabi-Yau varieties by the mirror theorem of \cite{Givental98} and \cite{LLY}. The periods for the LG models are related to relative Gromov-Witten invariants of the quasi-Fano varieties by the relative mirror theorem which is recently proved in \cite{FTY}.
\subsection{Definition}

In this section, we give a brief review of the definition of absolute and relative Gromov-Witten invariants.

Let $X$ be a smooth projective variety. We consider the moduli space $\overline{M}_{g,n,d}(X)$ of $n$-pointed, genus $g$, degree $d\in H_2(X,\mathbb Q)$ stable maps to $X$.
Let $\on{ev}_i$ be the $i$-th evaluation map, where
\begin{align*}
\on{ev}_i: \overline{M}_{g,n,d}(X) \rightarrow X, & \quad \text{for } 1\leq i\leq n.
\end{align*}

Let $s_i:\overline{M}_{g,n,d}(X)\rightarrow \mathcal C_{g,n,d}(X)$ be the $i$-th section of the universal curve, and let $\psi_i=c_1(s_i^*(\omega_{\mathcal C_{g,n,d}(X)/\overline{M}_{g,n,d}(X)}))$ be the descendant class at the $i$-th marked point.
Consider
\begin{itemize}
\item $\gamma_{i}\in H^*(X,\mathbb Q)$, for $1\leq i\leq n$;
\item  $a_{i}\in \mathbb Z_{\geq 0}$, for $1\leq i\leq n$.
\end{itemize}
Absolute Gromov-Witten invariants of $X$ are defined as
\begin{align}
\left\langle \prod_{i=1}^n \tau_{a_{i}}(\gamma_{i})\right\rangle^{X}_{g,n,d}:=\int_{[\overline{M}_{g,n,d}(X)]^{\on{vir}}}\psi_{1}^{a_{1}}\on{ev}^{*}_{1}(\gamma_{1})\cdots\psi_{n}^{a_{n}}\on{ev}^{*}_{n}(\gamma_{n}).
\end{align}
In this paper, we will be focusing on genus zero absolute Gromov-Witten invariants of Calabi-Yau varieties.

Let $X$ be a smooth projective variety and $D$ be a smooth divisor. We can study the relative Gromov-Witten invariants of $(X,D)$. 

For $d\in H_2(X,\mathbb Q)$, we consider a partition $\vec k=(k_1,\ldots,k_{n_1})\in (\mathbb Z_{>0})^{n_1}$ of $\int_d[D]$. That is,
\[
\sum_{i=1}^{n_1} k_i=\int_d[D].
\] 

We consider the moduli space $\overline{M}_{g,\vec k,n_2,d}(X,D)$ of $(n_1+n_2)$-pointed, genus $g$, degree $d\in H_2(X,\mathbb Q)$, relative stable maps to $(X,D)$ such that the contact orders for relative markings are given by the partition $\vec k$. We assume the first $n_1$ marked points are relative marked points and the last $n_2$ marked points are non-relative marked points. Let $\on{ev}_i$ be the $i$-th evaluation map, where
\begin{align*}
\on{ev}_i: \overline{M}_{g,\vec k,n_2,d}(X,D) \rightarrow D, & \quad\text{for } 1\leq i\leq n_1;\\
\on{ev}_i: \overline{M}_{g,\vec k,n_2,d}(X,D) \rightarrow X, & \quad \text{for } n_1+1\leq i\leq n_1+n_2.
\end{align*}

Write $\bar\psi_i= s^*\psi_i$ which is the class pullback from the corresponding descendant class on the moduli space $\overline{M}_{g,n_1+n_2,d}(X)$ of stable maps to $X$.
Consider
\begin{itemize}
\item $\delta_{i}\in H^*(D,\mathbb Q)$, for $1\leq i\leq n_1$;
\item $\gamma_{n_1+i}\in H^*(X,\mathbb Q)$, for $1\leq i\leq n_2$;
\item  $a_{i}\in \mathbb Z_{\geq 0}$, for $1\leq i\leq n_1+n_2$.
\end{itemize}
Relative Gromov-Witten invariants of $(X,D)$ are defined as
\begin{align}\label{relative-invariant-higher-dimension}
&\left\langle \prod_{i=1}^{n_1}\tau_{a_{i}}(\delta_{i})\left|\prod_{i=1}^{n_2} \tau_{a_{n_1+i}}(\gamma_{n_1+i})\right.\right\rangle^{(X,D)}_{g,\vec k,n_2,d}:=\\
\notag &\int_{[\overline{M}_{g,\vec k,n_2,d}(X,D)]^{\on{vir}}}\bar{\psi}_1^{a_1}\on{ev}^*_{1}(\delta_{1})\cdots \bar{\psi}_{n_1}^{a_{n_1}}\on{ev}^*_{n_1}(\delta_{n_1})\bar{\psi}_{n_1+1}^{a_{n_1+1}}\on{ev}^{*}_{n_1+1}(\gamma_{n_1+1})\cdots\bar{\psi}_{n_1+n_2}^{a_{n_1+n_2}}\on{ev}^{*}_{n_1+n_2}(\gamma_{n_1+n_2}).
\end{align}
We refer to \cite{IP}, \cite{LR}, \cite{Li01} and \cite{Li02} for more details about the construction of relative Gromov-Witten theory. Relative Gromov-Witten theory has recently been generalized to include negative contact orders (allowing $k_i<0$) in \cite{FWY} and \cite{FWY19}. We refer to \cite{FWY} and \cite{FWY19} for the precise definition of relative Gromov-Witten invariants with negative contact orders. In this paper, we will be focusing on genus zero relative Gromov-Witten invariants of quasi-Fano varieties along with their anticanonical divisors.

\subsection{A mirror theorem for Calabi-Yau varieties}\label{sec:mirror-theorem-CY}

In this section, we briefly review the mirror theorem for toric complete intersections following \cite{Givental98}. We are focusing on the case of Calabi-Yau varieties. A mirror theorem relates a generating function of Gromov-Witten invariants of a Calabi-Yau variety and periods of the mirror. The generating function considered by Givental \cite{Givental98} is called $J$-function. The periods are encoded in the so-called $I$-function in \cite{Givental98}.

The $J$-function is a cohomological valued function defined by
\[
J^X(\tau,z)=e^{\tau/z}\left(1+\sum_{\alpha}\sum_{d\in \on{NE}(X)_{\mathbb Z}\setminus \{0\}}Q^d\left\langle \frac{\phi_\alpha}{z(z-\psi)} \right\rangle_{0,1,d}^X\phi^\alpha \right),
\]
where the notation is explained as follows:
\begin{itemize}
\item $\tau=\sum_{i=1}^rp_i\log Q_i\in H^2(X)$.
\item $p_1,\ldots,p_r$ is an integral, nef basis of $H^2(X)$.
\item $\on{NE}(X)\subset H_2(X,\mathbb R)$ is the cone generated by effective curves and $\on{NE}(X)_{\mathbb Z}:=\on{NE}(X)\cap (H_2(X,\mathbb Z)/\on{tors})$.
\item $\{\phi_\alpha\}$ and $\{\phi^\alpha\}$ are dual bases of $H^*(X)$.
\item $\langle \frac{\phi_\alpha}{z(z-\psi)} \rangle_{0,1,d}^X$ is Gromov-Witten invariant of $X$ with degree $d$, 1-marked point  and the insertion is read as follows:
\[
\frac{\phi_\alpha}{z(z-\psi)}=\sum_{k\geq 0} \phi_\alpha\psi^k z^{-(k+2)}.
\]
\end{itemize} 

Let $X$ be a Calabi-Yau complete intersection in a toric variety $Y$ defined by a section of $E=L_0\oplus L_1\oplus\cdots \oplus  L_s$, where each $L_l$ is a nef line bundle. Recall that $\rho_l=c_1(L_l)$ and $-K_Y=\sum_{l=0}^s c_1(L_l)$. Consider the fiberwise $\mathbb C^*$-action on the total space $E$. Following Coates-Givental \cite{CG}, we consider the universal family 
\begin{displaymath}
    \xymatrix{ \mathcal C_{0,n,d} \ar[r]^{\on{ev}}\ar[d]^{\pi} &  Y\\
  \overline{M}_{0,n,d}(Y)&}.
\end{displaymath}
Coates-Givental \cite{CG} define
\[
E_{0,n,d}:=R\pi_*\on{ev}^*E\in K_{\mathbb T}^0(\overline{M}_{0,n,d}(Y)).
\]
Let $e(\cdot)$ denote the $\mathbb C^*$-equivariant Euler class. 
The twisted $J$-function is
\[
J^{e,E}(\tau,z)=e^{\tau/z}\left(1+\sum_{\alpha}\sum_{d\in \on{NE}(Y)_{\mathbb Z}\setminus \{0\}}Q^d\left\langle \frac{\phi_\alpha}{z(z-\psi)} \right\rangle_{0,1,d}^{e,E}\phi^\alpha \right),
\]
where
\[
\left\langle \prod_{i=1}^n \tau_{a_i}(\gamma_i)\right\rangle_{0,n,d}^{(e, E)}:=\int_{[\overline{M}_{0,n,d}(Y )]^{\on{vir}}}\prod_{i=1}^n(\on{ev}^*_i\gamma_i)\bar{\psi}_i^{a_i}e(E_{0,n,d})
\]
is called a twisted Gromov-Witten invariant of $Y$. The twisted $J$-function admits a non-equivariant limit $J^{Y,X}$, which satisfies
\[
i_*J^X(i^*(\tau),z)=J^{Y,X}(\tau,z)\cup \prod_{l=0}^s \rho_l,
\]
where $i:X\hookrightarrow Y$ is the inclusion.

Let $D_1,\ldots, D_m$ be the classes in $H^2(Y)$ which are Poincar\'e dual to the toric divisors in $Y$. The $I$-function is cohomological valued function defined by
\[
I^{X}(q,z)=e^{\sum_{i=1}^r p_i\log q_i/z}\sum_{d\in \on{NE}(Y)_{\mathbb Z}}\prod_{j=1}^m \left(\frac{\prod_{k=-\infty}^0(D_j+kz)}{\prod_{k=-\infty}^{\langle D_j,d\rangle}(D_j+kz)}\right)\left(\prod_{l=0}^s \prod_{k=1}^{\langle \rho_l, d\rangle }(\rho_l+kz)\right)q^d.
\]
The $I$-function can be expanded as a power series in $1/z$:
\[
I^X(q,z)=f_0(q)+\frac{f_1(q)}{z}+O(z^{-1}).
\]
The mirror theorem is stated as follows.
\begin{thm}[\cite{Givental98}, \cite{LLY}]
The $I$-function and the $J$-function are equal after change of variables:
\[
J^{Y,X}(f_1(q)/f_0(q),z)=I^X(q,z)/f_0(q).
\]
\end{thm}
\begin{remark}
Note that $I^X(q,z)$ takes values in $H^*(Y)$. The restriction $i^*I^X(q,z)$ takes values in $H^*(X)$ and lies in Givental's Lagrangian cone of the Gromov-Witten theory of $X$. A more common notation for our $I^X(q,z)$ is $I^{Y,X}(q,z)$. We choose not to use this notation to avoid any possible confusion with the relative $I$-function that we consider in Section \ref{sec:mirror-theorem-pairs}. 
\end{remark}

\begin{remark}
By setting $z=1$ for the $I$-functions, we obtain the bases of solutions to the Picard-Fuchs equations considered in previous sections.
\end{remark}
We list the $I$-functions for some examples of Calabi-Yau varieties that we consider in this paper.
\begin{example}
Let $Q_5$ be a smooth quintic threefold in $\mathbb P^4$. The $I$-function is
\[
I^{Q_5}(q,z)=e^{H\on{log}q/z}\sum_{d\geq 0}\left(\frac{\prod_{k=1}^{5d} (5H+kz)}{\prod_{k=1}^d(H+kz)^{5}}\right)q^d.
\]
\end{example}

\begin{example}
Let $\tilde{Q}_5$ be a complete intersection of bidegrees $(4,1)$ and $(1,1)$ in $\mathbb P^4\times \mathbb P^1$. The $I$-function is
\[
I^{\tilde{Q}_5}(q_1,q_0,z)=q_1^{H/z}q_0^{P/z}\sum_{d_1,d_0\geq0}\left(\frac{\prod_{k=1}^{4d_1+d_0} (4H+P+kz)\prod_{k=1}^{d_1+d_0}(H+P+kz)}{\prod_{k=1}^{d_1}(H+kz)^{5}\prod_{k=1}^{d_0}(P+kz)^{2}}\right)q_1^{d_1}q_0^{d_0}.
\]
\end{example}

\subsection{A mirror theorem for log Calabi-Yau pairs}\label{sec:mirror-theorem-pairs}

A mirror theorem for a relative pair $(X,D)$ is proved in \cite{FTY}. In this section, we focus on the case when $X$ is a toric complete intersection and $D$ is a smooth anticanonical divisor of $X$. The relative mirror theorem in \cite{FTY} is again formulated in terms of the relation between $J$-functions and $I$-functions under Givental's formalism for relative Gromov-Witten theory which is recently developed in \cite{FWY}. 

Following \cite{FWY}, the ring of insertions for relative Gromov-Witten theory is defined to be
\[
\HH=\bigoplus\limits_{i\in\mathbb Z}\HH_i,
\]
where $\HH_0=H^*(X)$ and $\HH_i=H^*(D)$ if $i\in \mathbb Z \setminus \{0\}$.
For an element $\alpha\in \HH_i$, we write $[\alpha]_i$ for its embedding in $\HH$. 
The pairing on $\HH$ is defined as follows:
\begin{equation}\label{eqn:pairing}
\begin{split}
([\alpha]_i,[\beta]_j) = 
\begin{cases}
0, &\text{if } i+j\neq 0,\\
\int_X \alpha\cup\beta, &\text{if } i=j=0, \\
\int_D \alpha\cup\beta, &\text{if } i+j=0, i,j\neq 0.
\end{cases}
\end{split}
\end{equation}

Let $t_{0,2}\in H^2(X)$ and $[t_{tw}]\in \HH\setminus \HH_0$. Following \cite{FTY}, we define the $J$-function.
\begin{defn}
The $J$-function for relative Gromov-Witten invariants of $(X,D)$ is a $\HH$-valued function
\begin{align*}
&J^{(X,D)}([t],z):=\\
&e^{t_{0,2/z}}\left(1+\frac{[t_{tw}]}{z}+\sum_{(n,d)\neq (0,0)}\sum_{\alpha}\frac{Q^d}{n!}\left\langle \frac{[\phi_\alpha]}{z(z-\bar\psi)},[t_{tw}],\ldots,[t_{tw}]\right\rangle^{(X,D)}_{0,n+1, d}[\phi^\alpha]\right),
\end{align*}
where $[t]=[t_{0,2}]+[t_{tw}]\in \HH$; $t_{0,2}=\sum_{i=1}^rp_i\log Q_i\in H^2(X)$; $\{[\phi_\alpha]\}$ and $\{[\phi^\alpha]\}$ are dual bases of $\HH$ under the pairing (\ref{eqn:pairing}).
\end{defn}

For simplicity, we only write down the $I$-function for the log Calabi-Yau pair $(X,D)$, where $X$ is a complete intersection in the toric variety $Y$. 
We consider the extended $I$-function for relative Gromov-Witten theory. It may be considered as a "limit" of the extended $I$-function for root stacks $X_{D,r}$ for $r\rightarrow \infty$. This is because of the relation between genus zero relative and orbifold Gromov-Witten invariants proved in \cite{ACW}, \cite{TY18} and \cite{FWY}. By \cite{FTY}, root stacks are hypersurfaces in toric stack bundles with fibers being weighted projective lines $\mathbb P^1_{1,r}$. Let $\sum_{i> 0}^b x_i[\textbf 1]_i$ represent the extended data. In the language of orbifold Gromov-Witten theory, $x_i$ corresponds to the box element of age $i/r$ in the extended stacky fan of $\mathbb P^1_{1,r}$. Note that the age $i/r$ in orbifold Gromov-Witten theory corresponds to contact order $i$ in relative Gromov-Witten theory. The limit is taken by first choosing $k_i$, $d$ and $b$, then taking a sufficiently large $r$. Since it works for any $b>0$, we may formally take $b\rightarrow \infty$ after taking $r\rightarrow \infty$. 

\begin{definition}\label{def-relative-I-function}
The extended $I$-function for relative invariants is defined as follows
\[
I^{(X,D)}(q,x,t,z)=I_++I_-,
\]
where
\begin{align*}
I_+:=&e^{\sum_{i=1}^r p_i\log q_i/z}\sum_{\substack{d\in \overline{\on{NE}}(Y),k_i\geq 0\\ \sum ik_i<\langle D, d\rangle } }\left(\prod_{j=1}^m\frac{\prod_{k=-\infty}^0(D_j+kz)}{\prod_{k=-\infty}^{\langle D_j,d\rangle}(D_j+kz)}\right)
\notag \left(\prod_{l=0}^s \prod_{k=1}^{\langle \rho_l, d\rangle }(\rho_l+kz)\right)\\
&\frac{\prod_{0<k\leq \langle D, d\rangle }(D+kz)}{D+(\langle D, d\rangle -\sum ik_i)z}\frac{\prod x_i^{k_i}}{z^{\sum k_i}\prod (k_i!)}[{\mathbf 1}]_{-\langle D, d\rangle +\sum ik_i}q^d,
\end{align*}
and 
\begin{align*}
I_-:=&e^{\sum_{i=1}^r p_i\log q_i/z}\sum_{\substack{d\in \overline{\on{NE}}(Y),k_i\geq 0\\ \sum ik_i\geq \langle D, d\rangle } }\left(\prod_{j=1}^{m}\frac{\prod_{k=-\infty}^0(D_j+kz)}{\prod_{k=-\infty}^{\langle D_j,d\rangle}(D_j+kz)}\right)
\notag \left(\prod_{l=0}^s \prod_{k=1}^{\langle \rho_l, d\rangle }(\rho_l+kz)\right)\\
&\left(\prod_{0<k\leq \langle D, d\rangle }(D+kz)\right)\frac{\prod x_i^{k_i}}{z^{\sum k_i}\prod (k_i!)}[{\mathbf 1}]_{-\langle D, d\rangle +\sum ik_i}q^d.
\end{align*}
\end{definition}

\begin{theorem}[\cite{FTY}, Theorem 1.5]
The extended relative $I$-function $I^{(X,D)}$ lies in Givental’s Lagrangian cone for relative invariants which is
defined in \cite{FWY}*{Section 7.5}.
\end{theorem}

For the purpose of this paper, we take the part of the $J$-function that takes values in $H^*(X)$. 
Therefore, we consider the function of the form
\begin{align*}
&J^{(X,D)}_0([t],z):=\\
&e^{t_{0,2}/z}z\left(1+\sum_{(n,d)\neq (0,0)}\sum_{\alpha}\frac{Q^d}{n!}\left\langle \frac{[\phi_\alpha]_0}{z-\bar\psi},[t_{tw}],\ldots,[t_{tw}]\right\rangle^{(X,D)}_{0,n+1, d}[\phi^\alpha]_0\right),
\end{align*}
where $[\phi_\alpha]_0$ and $[\phi^\alpha]_0$ are dual bases of $H^*(X)$. The corresponding $I$-function, denoted by $I_0$, is the part of $I_-$ that takes values in $\HH_0=H^*(X)$. Note that the $I$-function that we consider in Definition \ref{def-relative-I-function} is actually a non-equivariant limit of the twisted $I$-function from the ambient toric variety $Y$ similar to the absolute case in Section \ref{sec:mirror-theorem-CY}. Therefore, $I_0$ takes values in $H^*(Y)$ and $i^*I_0$ takes values in $H^*(X)$ where $i:X\hookrightarrow Y$ is the inclusion.

\begin{defn}
Let $X$ be a complete intersection in a toric variety $Y$ defined by a section of $E=L_0\oplus L_1\oplus\cdots \oplus  L_s$, where each $L_l$ is a nef line bundle. Let $\rho_l=c_1(L_l)$. Assuming $D$ is nef, the $I$-function is the $H^*(Y)$-valued function defined as
\begin{align}
I^{(X,D)}_0(q,x,z)=e^{\sum_{i=1}^r p_i\log q_i/z}\sum_{d\in \on{NE}(Y)_{\mathbb Z}}\sum_{k_i\geq 0, \sum ik_i=\langle D, d\rangle}\prod_{j=1}^m \left(\frac{\prod_{k=-\infty}^0(D_j+kz)}{\prod_{k=-\infty}^{\langle D_j,d\rangle}(D_j+kz)}\right)\\
\notag \left(\prod_{l=0}^s \prod_{k=1}^{\langle \rho_l, d\rangle }(\rho_l+kz)\right)\prod_{k=1}^{\langle D, d\rangle }(D+kz)\frac{\prod x_i^{k_i}}{z^{\sum k_i}\prod k_i!}q^d.
\end{align}
\end{defn}

\begin{rmk}
The nefness assumption on $D$ can be dropped when $D$ is coming from a toric divisor on $Y$. In this case, the relative Gromov-Witten theory of $(X,D)$ can be considered as a limit of the Gromov-Witten theory of the complete intersection on the root stack $Y_{D,r}$. Since $Y_{D,r}$ is a toric orbifold, the $I$-function for a complete intersection in toric orbifolds is known in \cite{CCIT14}.
\end{rmk}
Now, we consider some examples that we study in this paper.
\begin{example}
Let $X$ be a complete intersection of bidegrees $(4,1)$ in $\mathbb P^3\times \mathbb P^1$ and $D$ is its smooth anticanonical divisor. The $I$-function $I_0$ is
\begin{align*}
&I^{(X,D)}_0(q_1,q_0,x,z)=\\
&q_1^{H/z}q_0^{P/z}\sum_{d_1,d_0\geq0}\sum_{k_i\geq 0, \sum ik_i=d_0}\left(\frac{\prod_{k=1}^{4d_1+d_0} (4H+P+kz)}{\prod_{k=1}^{d_1}(H+kz)^{4}\prod_{k=1}^{d_0}(P+kz)}\right)\frac{\prod x_i^{k_i}}{z^{\sum k_i}\prod k_i!}q_1^{d_1}q_0^{d_0}.
\end{align*}
The extended $I$-function is
\[
I^{(X,D)}(q_1,q_0,x,z)=I_++I_-,
\]
where
\begin{align*}
I_+:=&q_1^{H/z}q_0^{P/z}\sum_{\substack{d_1,d_0\geq0,k_i\geq 0\\ \sum ik_i<d_0,}}\left(\frac{\prod_{k=1}^{4d_1+d_0} (4H+P+kz)}{\prod_{k=1}^{d_1}(H+kz)^{4}\prod_{k=1}^{d_0}(P+kz)}\right)\\
&\frac{1}{P+(d_0-\sum ik_i)z}\frac{\prod x_i^{k_i}}{z^{\sum k_i}\prod (k_i!)}[{\mathbf 1}]_{-d_0+\sum ik_i}q_1^{d_1}q_0^{d_0},
\end{align*}
and 
\begin{align*}
I_-:=&q_1^{H/z}q_0^{P/z}\sum_{\substack{d_1,d_0\geq0,k_i\geq 0 \\ \sum ik_i\geq d_0}}\left(\frac{\prod_{k=1}^{4d_1+d_0} (4H+P+kz)}{\prod_{k=1}^{d_1}(H+kz)^{4}\prod_{k=1}^{d_0}(P+kz)}\right)\\
&\frac{\prod x_i^{k_i}}{z^{\sum k_i}\prod (k_i!)}[{\mathbf 1}]_{-d_0+\sum ik_i}q_1^{d_1}q_0^{d_0}.
\end{align*}
The mirror map is given by the quotient of the coefficients of $z^0$ and $z^{-1}$. The $z^0$-coefficient is $\sum_{d_1\geq 0}\frac{(4d_1)!}{(d_1!)^4}q_1^{d_1}$. The coefficient of $z^{-1}$ is
\begin{align*}
&\sum_{d_1\geq 1}\frac{(4d_1)!}{(d_1!)^4}q_1^{d_1}\left(\sum_{k=d_1+1}^{4d_1}\frac 1k\right)(4H+P)+H\log q_1+P\log q_0\\
+&\sum_{d_1\geq 0,d_0>0}\frac{(4d_1+d_0)!}{(d_1!)^4(d_0)!d_0}[\mathbf 1]_{-d_0}q_1^{d_1}q_0^{d_0}+\sum_{i=1}^\rho\left(\sum_{d_1\geq 0, 0\leq d_0\leq i}\frac{(4d_1+d_0)!}{(d_1!)^4(d_0)!}q_1^{d_1}q_0^{d_0}[{\mathbf 1}]_{-d_0+i}\right)x_i.
\end{align*}

\end{example}

\begin{example}
Let $X$ be a complete intersection of bidegrees $(4,0)$ and $(1,1)$ in $\mathbb P^4\times \mathbb P^1$. Let $D$ be its smooth anticanonical divisor. The $I$-function $I_0$ is
\begin{align*}
&I^{(X,D)}_0(q_1,q_0,x,z)=\\
&q_1^{H/z}q_0^{P/z}\sum_{d_1,d_0\geq0}\sum_{k_i\geq 0, \sum ik_i=d_0}\left(\frac{\prod_{k=1}^{4d_1} (4H+kz)\prod_{k=1}^{d_1+d_0}(H+P+kz)}{\prod_{k=1}^{d_1}(H+kz)^{5}\prod_{k=1}^{d_0}(P+kz)}\right) \frac{\prod x_i^{k_i}}{z^{\sum k_i}\prod k_i!} q_1^{d_1}q_0^{d_0}.
\end{align*}
\end{example}

\begin{example}
Let $Q_4$ be the quartic threefold in $\mathbb P^4$ and $K3$ be the smooth anticanonical divisor. The $I$-function $I_0$ is
\[
I^{(Q_4,K3)}_0(q_1,x,z)=
e^{H\on{log}q_1/z}\sum_{d\geq 0}\sum_{k_i\geq 0, \sum ik_i=d}\left(\frac{\prod_{k=1}^{4d} (4H+kz)}{\prod_{k=1}^d(H+kz)^{4}}\right)\frac{\prod x_i^{k_i}}{z^{\sum k_i}\prod k_i!}q_1^d.
\]
\end{example}

\begin{example}
Let $X$ be the blow-up of $\mathbb P^3$ along the complete intersection of degrees $4$ and $5$ hypersurfaces. The extended $I$-function for relative invariants is defined as follows
\[
I^{(X,D)}(q_1,q_0,x,t,z)=I_++I_-,
\]
where
\begin{align*}
I_+:=&q_1^{H/z}q_0^{P/z}\sum_{\substack{d_1,d_0\geq0,k_i\geq 0\\ \sum ik_i<d_0-d_1,}}\left(\frac{\prod_{k=1}^{4d_1+d_0} (4H+P+kz)}{\prod_{k=1}^{d_1}(H+kz)^{4}\prod_{k=1}^{d_0}(P+kz)}\right)\\
&\frac{1}{P-H+(d_0-d_1-\sum ik_i)z}\frac{\prod x_i^{k_i}}{z^{\sum k_i}\prod (k_i!)}[{\mathbf 1}]_{-d_0+d_1+\sum ik_i}q_1^{d_1}q^{d_0},
\end{align*}
and 
\begin{align*}
I_-:=&q_1^{H/z}q_0^{P/z}\sum_{\substack{d_1,d_0\geq0,k_i\geq 0 \\ \sum ik_i\geq d_0-d_1}}\left(\frac{\prod_{k=1}^{4d_1+d_0} (4H+P+kz)}{\prod_{k=1}^{d_1}(H+kz)^{4}\prod_{k=1}^{d_0}(P+kz)}\right)\\
&\frac{\prod x_i^{k_i}}{z^{\sum k_i}\prod (k_i!)}[{\mathbf 1}]_{-d_0+d_1+\sum ik_i}q_1^{d_1}q_0^{d_0}.
\end{align*}
In particular, $I_0^{(X,D)}$ is the part of $I_-$ that takes values in $H^*(X)$ (that is, when $d_0-d_1=\sum ik_i$). We have
\begin{align*}
I_0^{(X,D)}(q_1,q_0,x,z):=q_1^{H/z}q_0^{P/z}\sum_{\substack{d_1,d_0\geq0,k_i\geq 0 \\ \sum ik_i= d_0-d_1}}\left(\frac{\prod_{k=1}^{4d_1+d_0} (4H+P+kz)}{\prod_{k=1}^{d_1}(H+kz)^{4}\prod_{k=1}^{d_0}(P+kz)}\right)
\frac{\prod x_i^{k_i}}{z^{\sum k_i}\prod (k_i!)}q_1^{d_1}q_0^{d_0}.
\end{align*}
\end{example}

\subsection{Relation with periods}
Since Gromov-Witten invariants are related to periods via mirror maps, the relation among periods given by the gluing formula implies a relation among different Gromov-Witten invariants via their respective mirror maps. 

Relative periods can be extracted from the $I$-functions for relative Gromov-Witten theory. We will illustrate this with an explicit example. Recall that, $\tilde{Q}_5$ is a complete intersection of bidegrees $(4,1)$ and $(1,1)$ in $\mathbb P^4\times \mathbb P^1$. The Calabi-Yau threefold $\tilde{Q}_5$ admits the following Tyurin degeneration
\[
\tilde{Q}_5\leadsto X_1\cup_{K3} X_2,
\]
where $X_1$ is a hypersurface of bidegree $(4,1)$ in $\mathbb P^3\times \mathbb P^1$ and $X_2$ is a complete intersection of bidegrees $(4,0), (1,1)$ in $\mathbb P^4\times \mathbb P^1$. Note that the $I$-function for relative invariants includes some extra variables $x_i$ corresponding to the contact order of relative markings. 

There are many ways to extract periods from the relative $I$-function. To simplify the discussion, we can extract the coefficient of $x_{d_0}q_1^{d_1}q_0^{d_0}$ for each $d_1,d_0\geq 0$ from $I_{0}^{(X_1,K3)}$ and $I_{0}^{(X_2,K3)}$. Setting $x_i=1$, then the resulting $I$-functions are
\begin{align*}
&I^{(X_1,K3)}_0(q_1,q_0,z)=\\
&q_1^{H/z}q_0^{P/z}\sum_{d_1,d_0\geq0}\left(\frac{\prod_{k=1}^{4d_1+d_0} (4H+P+kz)}{\prod_{k=1}^{d_1}(H+kz)^{4}\prod_{k=1}^{d_0}(P+kz)}\right)\frac{q_1^{d_1}q_0^{d_0}}{z^{\min(1,d_0)}},
\end{align*}
and 
\begin{align*}
&I^{(X_2,K3)}_0(q_1,q_0,z)=\\
&q_1^{H/z}q_0^{P/z}\sum_{d_1,d_0\geq0}\left(\frac{\prod_{k=1}^{4d_1} (4H+kz)\prod_{k=1}^{d_1+d_0}(H+P+kz)}{\prod_{k=1}^{d_1}(H+kz)^{5}\prod_{k=1}^{d_0}(P+kz)}\right)  \frac{q_1^{d_1}q_0^{d_0}}{z^{\min(1,d_0)}},
\end{align*}
which, after setting $z=1$, are exactly the bases of solutions for the Picard-Fuchs equations in Section \ref{sec:period-LG-example-1}. One can then extract the corresponding invariants from the $J$-functions with the corresponding mirror maps. The gluing formula for periods yields a formula relating the corresponding invariants by plugging in the corresponding $J$-functions along with mirror maps to Equation (\ref{Hadamard-I-function-conifold}).

\begin{remark}
The resulting gluing formula for Gromov-Witten invariants relates absolute Gromov-Witten invariants of the Calabi-Yau manifolds $X$ and $X_0$ to the relative Gromov-Witten invariants of the pairs $(X_1, X_0)$ and $(X_2, X_0)$.  This is clearly a sort of ``$B$-model motivated'' degeneration formula for genus zero Gromov-Witten invariants of Tyurin degenerations. Nevertheless, we do not know at present how to relate our gluing formula for periods to the degeneration formula for Gromov-Witten invariants as described in
\cites{Li01,Li02,LR}. 

More so than the complexity of the mirror maps and Birkhoff factorizations, there are several other obstructions to directly checking compatibility. For example, the relative mirror theorem only involves relative Gromov-Witten invariants whose insertions at the relative markings are cohomology classes pulled back from the ambient space. On the other hand, the degeneration formula involves relative Gromov-Witten invariants whose insertions at relative markings having cohomology classes that are not pulled back from the ambient space. Moreover, there may be relative invariants with negative contact orders involved in the relative mirror theorem, while the degeneration formula only has relative invariants with markings of positive contact orders.
\end{remark}
\bibliographystyle{amsxport}
\bibliography{universalreferences.bib}

\end{document}